\documentclass{amsart}
\usepackage{amssymb}
\usepackage{mathrsfs}
\usepackage{stmaryrd}
\usepackage{bbm}
\usepackage{oldgerm}
\usepackage[english]{babel}
\usepackage[T1]{fontenc}
\usepackage[latin1]{inputenc}
\usepackage[all]{xy}
\usepackage{hyperref}
\usepackage{upref}
\usepackage{xcolor}
\usepackage{tikz-cd}

\hypersetup{
 colorlinks,
 linkcolor={red!50!black},
 citecolor={blue!50!black},
 urlcolor={blue!80!black}
}

\newtheorem{teo}[subsection]{Theorem}
\newtheorem{prop}[subsection]{Proposition}
\newtheorem{cor}[subsection]{Corollary}
\newtheorem{lem}[subsection]{Lemma}

\theoremstyle{definition}

\newtheorem{defi}[subsection]{Definition}
\newtheorem{rema}[subsection]{Remark}

\newtheorem{exemple}[subsection]{Example}
\newtheorem{exemples}[subsection]{Examples}
\newtheorem{constr}[subsection]{Construction}

\numberwithin{equation}{subsection}

\newcommand{\ka}{\kappa}

\newcommand{\M}{\mathcal{M}}

\newcommand{\m}{\mathfrak{m}}
\newcommand{\Ow}{\mathcal{O}}
\newcommand{\p}{\mathfrak{p}}

\newcommand{\pur}{\mathrm{pur}}

\title{A new proof of Raynaud-Gruson's flattening theorem}
\author{Quentin Guignard}
\address{Institut des
Hautes \'Etudes Scientifiques, 35 route de Chartres, 91440 Bures-sur-Yvette, France}
\address{
\'Ecole Normale Sup\'erieure, 45 rue d'Ulm, 75005 Paris, France}
\email{quentin.guignard@ens.fr}

\begin{document}

\begin{abstract}
We give a new proof of Raynaud-Gruson's theorem regarding flattening by blow-ups. The proof is direct, by working directly on the inverse limit of admissible blow-ups, which is a valuative space similar to the classical Zariski-Riemann space. These valuative spaces are defined in a broader context and always have an analogous flattening property; this yields a generalization of Raynaud-Gruson's theorem. 
\end{abstract}

\maketitle
\tableofcontents

\section{Introduction \label{chap1introrgg}}

In Raynaud's approach to rigid geometry, the category of quasi-compact quasi-separated rigid spaces over a completely valued field $K$ of rank $1$ is defined as the localization with respect to ``admissible blow-ups'' of the category of finitely presented formal schemes over the ring of integers $R \subseteq K$ (cf. \cite{Raynaud}). It has been known since then that a flat morphism of quasi-compact quasi-separated rigid spaces over $K$ can be represented by a flat morphism between appropriate $R$-models. A precise statement and a proof can be found in (\cite{BL93}, 5.2) or (\cite{Abbes11}, 5.8.1). The schematic version of this flattening theorem had been previously proved by Raynaud and Gruson in 1971:

\begin{teo}[(\cite{Stacks} $\mathrm{0815}$), cp. (\cite{Raynaud-Gruson71} $I.5.2.2$)]\label{chap1RGtheo} Let $Y$ be a scheme of finite presentation over a quasi-compact and quasi-separated scheme $X$, let $U$ be a quasi-compact open subset of $X$ and let $\mathcal{F}$ be a quasi-coherent $\Ow_Y$-module of finite type. Assume that the restriction of $\mathcal{F}$ to $Y \times_X U$ is a finitely presented $\Ow_{Y \times_X U}$-module which is flat over $U$. Then there exists a blow-up $f : X' \rightarrow X$ such that:
\begin{itemize}
\item[$(1)$] The center of the blow-up $f$ is a finitely presented closed subscheme of $X$, disjoint from $U$.
\item[$(2)$] If $Y'$ is the strict transform of the $X$-scheme $Y$ along $f$, then the strict transform $\mathcal{F}'$ of $\mathcal{F}$ along $f$ is finitely presented over $\Ow_{Y'}$ and flat over $X'$.
\end{itemize}
\end{teo}

Here, by ``strict transform'' we mean the following: if $Z$ is the exceptional Cartier divisor of the blow-up $f$, then $Y'$ is the closed subscheme of $X' \times_X Y$ defined by the vanishing of the quasi-coherent ideal of sections supported on $Z \times_X Y$. Similarly, the strict transform $\mathcal{F}'$ of $\mathcal{F}$ along $f$ is the pullback to $Y'$ of the quotient of $\mathcal{F} \otimes_{\Ow_Y} \Ow_{X' \times_X Y}$ by the submodule of its sections supported on $Z \times_X Y$.

Our proof of Raynaud-Gruson's theorem proceeds in two steps.We first consider the (filtered) limit
$$\widetilde{X} = \lim_{X' \rightarrow X} \ X' $$
of all blow-ups as allowed in the statement of the theorem. This is not a scheme, but a ringed space which we identify to a valuative space \ref{chap1comparisoneclate}. These valuative spaces are constructed and studied in greater generality in Section \ref{chap1valspaces}, using the langage of ``rings with constructible supports'', whose set-up is the purpose of Sections \ref{chap1constrrings} and \ref{chap1phiringssheaves}. We obtain in this broader context a more general version of Raynaud-Gruson's theorem, namely Theorem \ref{chap1RG2}, which implies that a variant of the conclusion of Theorem \ref{chap1RGtheo} holds over $\widetilde{X}$ (cf. Section \ref{chap1flattprop}). The second and last step consists in descending the properties shown to hold on the limit $\widetilde{X}$ to some blow-up (cf. Section \ref{chap1proofofrg}). We now give a more detailed description of the content of this text.

In Section \ref{chap1constrrings}, we define the notion of ``rings with constructible supports'' or ``$\Phi$-rings'' (cf. \ref{chap11.2}) and we consider their basic properties. We include a few sorites on the notion of depth (cf. \ref{chap11.5} and \ref{chap11.6}), and give an interpretation of the closure and purification functors from (\cite{EGA42}, $5.9$, $5.10$) as adjoint functors, cf. \ref{chap11.9} and \ref{chap11.9bis}. We conclude Section \ref{chap1constrrings} with the definition and the study of the notion of ``$\Phi$-local'' $\Phi$-rings.

Section \ref{chap1phiringssheaves} is devoted to the definition of sheaves of $\Phi$-rings and to their basic properties. It is noticed there that the notion of $\Phi$-ring is well-defined in any topos, cf. \ref{chap10.2.0.1}, and that the notions and properties from Section \ref{chap1constrrings} extend to this broader context.

In Section \ref{chap1valspaces}, to any locally $\Phi$-ringed topological space $X$ we associate a valuative space $\widetilde{X}$, which is itself a locally $\Phi$-ringed topological space, endowed with a morphism $\widetilde{X} \rightarrow X$, cf. \ref{chap10.2.1}, \ref{chap10.2.2} and \ref{chap10.2.3}. We call this valuative space the ``$\Phi$-localization'' of $X$. It is actually a $\Phi$-locally $\Phi$-ringed topological space, and the morphism $\widetilde{X} \rightarrow X$ is universal for this property, cf. \ref{chap10.2.4}. We then consider the context of Theorem \ref{chap1RGtheo}, in which we endow the base scheme $X$ with a structure of locally $\Phi$-ringed topological space; we show that the valuative space associated to $X$ coincides with the limit of all blow-ups as allowed in Theorem \ref{chap1RGtheo}, cf. \ref{chap1comparisoneclate}.

In Section \ref{chap1flattprop}, $\Phi$-localizations are shown to have a flattening property, cf. \ref{chap1RG1} and \ref{chap1RG2}. The main argument is purely local, cf. \ref{chap10.3.2}, and rests upon the equational criterion for flatness. Theorem \ref{chap1RGtheo} is then deduced from \ref{chap1RG2} in Section \ref{chap1proofofrg}.

%In the first step, the limit over all blow-ups allowed in the theorem is not a scheme, but a valuative space together with a structural sheaf which are constructed and studied in Section \ref{chap1valspaces}. These valuative spaces naturally appear in the broader context of ``rings with constructible supports'', whose set-up is the purpose of Sections \ref{chap1constrrings} and \ref{chap1phiringssheaves}. This ultimately results in a more general version of Raynaud-Gruson's theorem, namely Theorem \ref{chap1RG2}. 

Ahmed Abbes has brought to our attention during the preparation of this text that a proof of Raynaud-Gruson's theorem similar to ours has been announced by Kazuhiro Fujiwara and Fumiharu Kato in \cite{FujiKato}. Besides the sketch given in (\cite{FujiKato} 5.6), their proof does not seem to have appeared in print. We notice that our treatment of the local case, namely \ref{chap10.3.2} and \ref{chap12.2.5}, seems to differ from theirs.

\subsection*{Acknowledgements} This work is part of the author's PhD dissertation and it was prepared at the Institut des Hautes \'Etudes Scientifiques and the \'Ecole Normale Sup\'erieure while the author benefited from their hospitality and support. The author is indebted to Ahmed Abbes and the anonymous referee for their thorough review of this text and for their numerous suggestions and remarks. Further thanks go to Lucia Mocz, Christophe Soul\'e and Salim Tayou for their review of an earlier version of this text, and to Fu Lei and Fabrice Orgogozo for several corrections and remarks.
 
\section{Rings with constructible supports \label{chap1constrrings}}

\subsection{\label{chap111}}  Let $A$ be a ring. A \textbf{family of constructible supports} on $A$ is a set $\Phi$ of finitely generated ideals of $A$ which satisfies the following two properties:
\begin{itemize}
\item[$(1)$] $\Phi$ is stable under finite products of ideals. In particular the unit ideal belongs to $\Phi$.
\item[$(2)$] If a finitely generated ideal $I$ of $A$ contains an element of $\Phi$, then $I$ is itself in $\Phi$.
\end{itemize}

\begin{exemples} 

\begin{itemize}\label{chap1exempleadique}

\item[$(i)$] Given an arbitrary set $\Phi_0$ of ideals of $A$, the \textbf{family of constructible supports generated by $\Phi_0$} is the set of finitely generated ideals of $A$ which contain a finite product of elements of $\Phi_0$. If $\Phi_0$ consists only of finitely generated ideals, then this family of constructible supports is the smallest one containing $\Phi_0$. 
\item[$(ii)$] If $A$ is a preadmissible topological ring in the sense of (\cite{EGA1} $0.7.1.2$), i.e. a linearly topologized ring with an ideal of definition (an open ideal $I$ of $A$ such that any neighbourhood of $0$ in $A$ contains $I^N$ for some integer $N$), then one can consider the family of constructible supports generated by all the ideals of definition. This family of constructible supports consists precisely of the finitely generated open ideals of $A$ whenever $A$ is preadic in the sense of (\cite{EGA1} $0.7.1.9$), i.e. when $A$ has an ideal of definition $I$ such that $I^N$ is open in $A$ for any integer $N$.
\end{itemize}
\end{exemples}

\subsection{\label{chap11.2}} The category $\Phi \mathrm{Rings}$ of \textbf{rings with constructible supports} or \textbf{$\Phi$-rings} is defined as follows:
\begin{itemize} 
\item[$\triangleright$] Its objects are the pairs $(A,\Phi_A)$ consisting of a ring $A$ and a family of constructible supports $\Phi_A$ on $A$ (cf. \ref{chap111}). An ideal of $A$ is said to be \textbf{admissible} if it belongs to $\Phi_A$.
\item[$\triangleright$] Its morphisms are the ring homomorphisms $f: A \rightarrow B$ such that for any admissible ideal $I$ of $A$, the ideal $f(I)B$ of $B$ is admissible.
\end{itemize}

From now on, we commit the abuse of notation where we use the same name for a ring with constructible supports and for its underlying ring, just as it is customary to use the same symbol for a ring and for its underlying set. 

\begin{prop}\label{chap1compcocomp} The category $\Phi \mathrm{Rings}$ is complete and cocomplete, and the forgetful functor from $\Phi \mathrm{Rings}$ to the category of rings preserve all small limits and colimits.
\end{prop}

Indeed let us consider a functor $F$ from a small category $J$ to $\Phi \mathrm{Rings}$. Let us endow the ring $\lim F$ with the family of constructible supports which consists of all the finitely generated ideals $I$ such that $I F(j)$ is admissible in $F(j)$ for all $j$. The resulting ring with supports is a limit of $F$ in $\Phi \mathrm{Rings}$.
Moreover if $\mathrm{colim}\ F$ is endowed with the family of constructible supports generated in the sense of \ref{chap1exempleadique} $(i)$ by all the ideals generated by an admissible ideal of $F(j)$ for some $j$, then the resulting ring with supports is a colimit of $F$ in $\Phi \mathrm{Rings}$.

\begin{prop}\label{chap1commute} The forgetful functor from $\Phi \mathrm{Rings}$ to the category of rings has both a left adjoint and a right adjoint.
\end{prop}

Indeed the functor which endows a ring $A$ with the family of constructible supports consisting only of the unit ideal is a left adjoint to the forgetful functor, while the functor which endows a ring $A$ with the family of constructible supports consisting of all of its finitely generated ideals is a right adjoint to the forgetful functor.

\begin{rema} \label{chap1exempleadique2}
Let $f: A \rightarrow B$ be a ring homomorphism between preadic topological rings in the sense of (\cite{EGA1} $0.7.1.9$). Let us endow $A$ and $B$ with the family of constructible supports consisting of all their finitely generated open ideals as in Example \ref{chap1exempleadique} $(ii)$, and let us assume that $A$ has a finitely generated ideal of definition. Then $f$ yields a morphism of $\Phi \mathrm{Rings}$ if and only if the topology on $B$ is finer that the topology induced by that of $A$. In particular $f$ may be a morphism of $\Phi \mathrm{Rings}$ without being continuous.
\end{rema}

\subsection{\label{chap11.3}} A \textbf{family of supports} on a topological space $X$ is a set $\Phi$ of closed subsets of $X$ which satisfies the following two properties:
\begin{itemize}
\item[$(1)$] $\Phi$ is stable under finite unions. In particular the empty set belongs to $\Phi$.
\item[$(2)$] If a closed subset of $Z$ of $X$ is contained in an element of $\Phi$, then $Z$ is itself in $\Phi$.
\end{itemize}
This definition is taken from (\cite{Godement73} II.2.5), besides the fact that we require a family of supports to be nonempty.\footnote{It is indeed mistakenly asserted in (\cite{Godement73} II.2.5) that for a family of supports $\Phi$ in the sense defined there and for a sheaf of abelian groups on $X$, the set of global sections of this sheaf with supports in $\Phi$ is an abelian group. This requires $\Phi$ to contain the support of the zero section, hence to be nonempty.} The notion of family of supports was introduced in \cite{Cartan50}, but additional conditions on $\Phi$ were required, e.g. paracompactness of its elements.

\begin{exemple} \label{chap1exemplesupport}
Given an arbitrary set $\Phi_0$ of closed subsets of $X$, the set of closed subsets of $X$ which are contained in a finite union of elements of $\Phi_0$ is a family of supports on $X$. It is the \textbf{family of supports generated by $\Phi_0$}, i.e. the smallest family of supports on $X$ containing $\Phi_0$
\end{exemple}

\subsection{\label{chap11.3.1}} If $\Phi$ is a family of constructible supports on a ring $A$ in the sense of \ref{chap11.2}, then the set of closed subsets $Z \subseteq X$ contained in $V(I)$ for some $I \in \Phi$ is a family of supports on $\mathrm{Spec}(A)$ which is generated (cf. \ref{chap1exemplesupport}) by its globally constructible elements (\cite{EGA1}, 0.2.3.2). Conversely, if $\Phi$ is a family of supports on $\mathrm{Spec}(A)$ generated by its globally constructible elements then the set of finitely generated ideals $I$ such that $V(I)$ belongs to $\Phi$ is a family of constructible supports on $A$.

These two constructions are inverse to each other. Indeed if $J$ is a finitely generated ideal such that $V(J) \subseteq V(I)$ for some admissible ideal $I =(f_1,\dots,f_r)$, then $f_i^{N_i}$ belongs to $J$ for some integer $N_i$. For $N = \sum_i N_i$ this yields $I^N \subseteq J$, so that $J$ is admissible. On the other hand, if we start from a family of supports $\Phi$ on $\mathrm{Spec}(A)$ generated by its globally constructible elements, then for any globally constructible element $Z$ of $\Phi$, the complement of $Z$ is a globally constructible open subset of $\mathrm{Spec}(A)$, hence is quasi-compact, so that it is a finite union of standard open subsets $(D(f_i))_{1 \leq i \leq r}$. Thus $Z = V(I)$ where $I$ is the ideal generated by $(f_i)_{1 \leq i \leq r}$.

\subsection{\label{chap11.4}} Let $X$ be a topological space and let $\Phi$ be a family of supports (cf. \ref{chap11.3}). Let us consider as in (\cite{Godement73} $II.2.5$) the functor $\Gamma_{\Phi}(X,-)$ which to a sheaf of abelian groups $\mathcal{F}$ on $X$ associates the abelian group
$$
\Gamma_{\Phi}(X,\mathcal{F}) = \{ s \in \Gamma(X,\mathcal{F}) \ | \ \mathrm{supp}(s) \in \Phi \}.
$$
It is a left exact additive functor between abelian categories, and its source has enough injectives (\cite{Godement73} $7.1.1$). Consequently, it has right derived functors (\cite{Godement73} 7.2), which are denoted by $H_{\Phi}^{q}(X,\mathcal{F})$. When $\Phi$ is the set of closed subsets of a given closed set $Z$, these groups are also denoted by $H_{Z}^{q}(X,\mathcal{F})$.

\begin{prop}\label{chap1comp} For any sheaf of abelian groups $\mathcal{F}$ on $X$, one has
$$
H_{\Phi}^{q}(X,\mathcal{F}) \cong \mathrm{colim} \ H_{Z}^{q}(X,\mathcal{F}),
$$
where the (filtered) colimit runs over elements $Z$ of $\Phi$.
\end{prop}

Indeed, if $\mathcal{F} \rightarrow \mathcal{J}^{\bullet}$ is an injective resolution of $\mathcal{F}$, then $\Gamma_{\Phi}(X,\mathcal{J}^{\bullet})$ is the (cofiltered) colimit of the complexes of abelian groups $\Gamma_{Z}(X,\mathcal{J}^{\bullet})$ when $Z$ runs over elements of $\Phi$. Since filtered colimits are exact, taking cohomology groups of these complexes yields \ref{chap1comp}.

 \begin{prop}\label{chap10.2top} Let $Z$ and $Z'$ be closed subsets of $X$, and let $\mathcal{F}$ be a sheaf of abelian groups on $X$.
\begin{itemize}
\item[$(a)$] There exists a long exact sequence
$$
\cdots \rightarrow H_{Z \cap Z'}^q(X,\mathcal{F}) \rightarrow H_Z^q(X,\mathcal{F}) \oplus H_{Z'}^q(X,\mathcal{F}) \rightarrow H_{Z \cup Z'}^q(X,\mathcal{F}) \rightarrow H_{Z \cap Z'}^{q+1}(X,\mathcal{F}) \rightarrow \cdots
$$
\item[$(b)$] Let $U = X \setminus Z$. There exists a long exact sequence
$$
\cdots \rightarrow H_{Z \cap Z'}^q(X,\mathcal{F}) \rightarrow H_{Z'}^q(X,\mathcal{F}) \rightarrow H_{Z' \cap U}^q \left( U, \mathcal{F}_{| U} \right) \rightarrow H_{Z \cap Z'}^{q+1}(X,\mathcal{F}) \rightarrow \cdots
$$
\end{itemize}
\end{prop}

Let $\mathcal{F} \rightarrow \mathcal{J}^{\bullet}$ be a flasque resolution (\cite{Kashiwara-Shapira94} II.$2.4.6vi$). Part $(a)$ follows from the exactness of the sequence of complexes
$$
0 \rightarrow \Gamma_{Z \cap Z'}(X,\mathcal{J}^{\bullet}) \rightarrow \Gamma_{ Z}(X,\mathcal{J}^{\bullet}) \oplus  \Gamma_{Z'}(X,\mathcal{J}^{\bullet}) \rightarrow \Gamma_{Z \cup Z'}(X,\mathcal{J}^{\bullet}) \rightarrow 0,
$$
while part $(b)$ follows from the exactness of the sequence of complexes
$$
0 \rightarrow \Gamma_{Z \cap Z'}(X,\mathcal{J}^{\bullet}) \rightarrow \Gamma_{Z'}(X,\mathcal{J}^{\bullet}) \rightarrow \Gamma_{Z' \cap U}(U,\mathcal{J}^{\bullet}_{|U}) \rightarrow 0.
$$ 
Alternatively, part $(b)$ is an instance of (\cite{SGA2} $2.2.8$, $2.2.2$).

\subsection{\label{chap11.5}} Let $A$ be a $\Phi$-ring (cf. \ref{chap11.2}), and let $\Phi_A$ be the corresponding family of supports on $X = \mathrm{Spec}(A)$ as in \ref{chap11.3}. If $M$ is an $A$-module, corresponding to the quasi-coherent sheaf $\widetilde{M}$ on $X$, then we set
$$
H_I^q(M)  =  H^q_{V(I)}(X,\widetilde{M}) \ \ \mathrm{ and } \ \ H_{\Phi_A}^q(M) = H_{\Phi_A}^q(X,\widetilde{M}).
$$ 

\begin{defi}\label{chap10.5} Let $d \geq 0$ be an integer. A module $M$ over the $\Phi$-ring $A$ is said to be $d$-deep if $H_{\Phi_A}^q(M)$ vanishes for each integer $q < d$. 
\end{defi}

\begin{rema} Assume that $A$ is noetherian, that its family of constructible supports is generated by a single ideal $I$, and that $M$ is finitely generated. The $I$-depth of $M$ in the sense of (\cite{SGA2} III.2.3) is at least $d$ if and only if $M$ is $d$-deep in the sense of the definition \ref{chap10.5}. This follows from (\cite{SGA2} III.3.1, III.3.3).
\end{rema}

\begin{lem}\label{chap10.5.1} Let $d \geq 0$ be an integer. Let $I$ be a finitely generated ideal of a ring $A$, and let $M$ be an $A$-module such that $H^q_I(M) = 0$ for any $q < d$. Then for any finitely generated ideal $J$ containing $I$, we have $H^q_J(M) = 0$ for any $q < d$, and the canonical morphism $H^d_J(M) \rightarrow H^d_I(M)$ is injective.
\end{lem}

By induction on the number of generators of $J$, one can assume that $J = (I,g)$ for some element $g$. One then applies \ref{chap10.2top} $(b)$ to $Z = V(g)$ and $Z' = V(I)$, so that $U = D(g)$. The resulting long exact sequence takes the form
$$
\cdots \rightarrow H_{J}^q(M) \rightarrow H_I^q(M) \rightarrow H_{I}^q \left( M [ g^{-1} ] \right) \rightarrow H_{J}^{q+1}(M) \rightarrow \cdots.
$$
By (\cite{SGA2} II.2) and by the hypothesis, one has 
$$
H_{I}^q \left( M [ g^{-1} ]  \right)  \cong H_{I}^q \left( M \right) [ g^{-1} ]  \cong 0
$$
for $q < d$. This proves the injectivity of $H_J^{q}(M) \rightarrow H_I^{q}(M)$ for $q \leq d$ and concludes the proof.

\begin{prop}\label{chap10.5bis} Let $A$ be a $\Phi$-ring and let $M$ be a $d$-deep $A$-module. Then $H_I^q(M)=0$ for any admissible ideal $I$ and any $q < d$.
\end{prop}

We prove Proposition \ref{chap10.5bis} by induction on $d$, the case $d=0$ being tautological. We are thus led to assume that $d \geq 1$ and that the result has been proved for $(d-1)$-deep modules. By \ref{chap1comp} one has 
$$
H_{\Phi_A}^q(M) \cong \mathrm{colim} \ H_I^q(M),
$$
where $I$ runs over the admissible ideals of $A$. It is therefore sufficient to prove that for any pair $I \subseteq J$ of admissible ideals and for any $q<d$ the homomorphism $H_J^q(M) \rightarrow H_I^q(M)$ is injective. The latter fact follows from Lemma \ref{chap10.5.1}, since $H_I^q(M) = 0$ for $q < d-1$ and for any admissible ideal $I$ of $A$.

%By induction on the number of generators of $J$, one can assume that $J = (I,g)$ for some element $g$. One then applies \ref{chap10.2top} $(b)$ to $Z = V(g)$ and $Z' = V(I)$, so that $U = D(g)$. The resulting long exact sequence takes the form
%$$
%\cdots \rightarrow H_{J}^q(M) \rightarrow H_I^q(M) \rightarrow H_{I}^q \left( M [ g^{-1} ] \right) \rightarrow H_{J}^{q+1}(M) \rightarrow \cdots.
%$$
%By (\cite{SGA2} II.2) and the induction hypothesis, one has 
%$$
%H_{I}^q \left( M [ g^{-1} ]  \right)  \cong H_{I}^q \left( M \right) [ g^{-1} ]  \cong 0
%$$
%for $q < d-1$. This proves the injectivity of $H_J^{q}(M) \rightarrow H_I^{q}(M)$ for $q <d$ and concludes the proof.

\begin{cor}\label{chap10.5ter} Let $A$ be a $\Phi$-ring and let $\Phi_0$ be a set of finitely generated ideals of $A$ generating the family of constructible supports of $A$ as in \ref{chap1exempleadique}$(i)$. An $A$-module $M$ is $d$-deep if and only if $H_I^q(M)=0$ for any $I$ in $\Phi_0$ and any $q < d$.
\end{cor}

The forward direction follows from \ref{chap10.5bis}. For the converse we proceed by induction on $d$, the case $d=0$ being tautological. We first note that $H_I^q(M)=0$ for any $q < d$ whenever $I$ is a finite product of elements of $\Phi_0$. Indeed, if $I$ and $I'$ satisfy $H_I^q(M) = H_{I'}^q(M)=0$ for any $q < d$, then by \ref{chap10.2top}$(a)$ we have an exact sequence
$$
0 \rightarrow H^q_{I I'}(M) \rightarrow  H^{q+1}_{I +I'}(M) \rightarrow H^{q+1}_{I }(M) \oplus H^{q+1}_{I'}(M),
$$
for any $q < d$. Since the canonical homomorphism $H^{q+1}_{I +I'}(M) \rightarrow H^{q+1}_{I }(M) $ is injective by Lemma \ref{chap10.5.1}, we obtain the vanishing of $H^q_{I I'}(M)$ for any $q <d$.

 Any admissible ideal $J$ must contain an ideal $I$ which is a finite product of elements of $\Phi_0$. By Lemma \ref{chap10.5.1}, we obtain $H_J^q(M)=0$ for any admissible ideal $J$ and any $q < d$. By \ref{chap1comp} the $A$-module $M$ must be $d$-deep.

\begin{prop}\label{chap1interest} Let $M$ be a $1$-deep $A$-module, and let $I = (f_1,\dots,f_r)$ be an admissible ideal of $A$. Then the homomorphism
$$
\mathrm{Hom}_A(I,M) \rightarrow M^{\oplus r}
$$
sending an element $\psi$ to $(\psi(f_i))_{1 \leq i \leq r}$ is injective, and its image consists of the $r$-uples $(m_i)_{1 \leq i \leq r}$ such that $f_i m_j = f_j m_i$ for all $i,j$.
\end{prop}

As $I$ is generated by $f_1,\dots,f_r,$ the injectivity is clear. In order to characterize the image, let us consider a free $A$-module $F$ of rank $r$ with basis $(e_i)_{1 \leq i \leq r}$ and let $G$ be its quotient by the relations $f_i e_j - f_j e_i$. Consider the homomorphism $G \rightarrow I$ sending $e_i$ to $f_i$, and let $H$ be its kernel. The exact sequence
$$
0 \rightarrow H \rightarrow G \rightarrow I
$$
yields the exact sequence
$$
\mathrm{Hom}_A(I,M) \rightarrow \mathrm{Hom}_A(G,M) \rightarrow \mathrm{Hom}_A(H,M) \rightarrow 0.
$$
We thus have to show that $\mathrm{Hom}_A(H,M)$ vanishes. Since $M$ has no nonzero $I$-torsion, it is sufficient to show that $IH = 0$. But if $x = \sum_i a_i e_i$ belongs to $H$ then
$$
f_j x = f_j \left( \sum_i a_i e_i\right) = \left(\sum_i a_i f_i\right)  e_j = 0
$$
holds in $G$ for any $j$, so that $Ix = 0$.

\begin{prop}\label{chap10.1} Let $M$ be an $A$-module. If $I$ is an ideal of $A$ generated by $f_1,\dots,f_r$, then $H_I^q(M)$ is isomorphic to the $q$-th cohomology group of the \v{C}ech complex
$$
C^{\bullet}(M,f_{\bullet}) \ : \ 0 \rightarrow M \rightarrow \prod_i M\left[ \frac{1}{f_i} \right] \rightarrow \prod_{i<j} M \left[ \frac{1}{f_i f_j} \right] \rightarrow \cdots \rightarrow M\left[ \frac{1}{f_1 \dots f_r} \right] \rightarrow 0
$$
where $M$ is in degree $0$. This isomorphism is functorial in $M$.
\end{prop}

This is (\cite{SGA2} II.5).

\begin{rema} It follows from Proposition \ref{chap10.1} that for any finitely generated ideal $I$ of $A$, any $A$-module $M$, and any flat ring homomorphism $A \rightarrow A'$, the base change homomorphism
$$
H_I^q(M) \otimes_A A' \rightarrow H_{I A'}^q(M \otimes_A A')
$$ 
is bijective for each integer $q$. 
\end{rema}

\begin{cor}\label{chap1transitive} Let $f : A \rightarrow B$ be a homomorphism of $\Phi$-rings, an let $M$ be a $B$-module. If $M$ is $1$-deep (respectively, $2$-deep) as a $B$-module, then it is so as an $A$-module. The converse holds if the family of constructible supports of $B$ is generated by that of $A$.
\end{cor}

Indeed, it follows from Proposition \ref{chap10.1} that for any admissible ideal $I$ of $A$ and any integer $q$, we have 
$$
H^q_I(M) \cong H^q_{IB}(M),
$$
where $M$ is considered as an $A$-module on the left hand side, and as a $B$-module on the right hand side.

\begin{cor}\label{chap11deep} An $A$-module $M$ is $1$-deep if and only if for any admissible ideal $I$ of $A$, $M$ has no nonzero $I$-torsion.
\end{cor}

Indeed by Proposition \ref{chap10.1} (or by definition) the module $H_I^0(M)$ consists of all $m$ in $M$ such that $f_i^{N_i} m = 0$ for some integers $(N_i)_{1 \leq i \leq r}$. For $N = \sum_i N_i$ this gives $I^N m =0$, so that $H_{\Phi_A}^0(M)$ is the submodule of elements $m$ in $M$ such that $Im =0$ for some admissible ideal $I$.

\begin{cor}\label{chap12deep} An $A$-module $M$ is $2$-deep if and only if for any admissible ideal $I$ of $A$, the homomorphism
$$
M \rightarrow \mathrm{Hom}_A(I,M)
$$
which sends an element $m$ to $(x \mapsto xm)$ is bijective.
\end{cor}
 
Assume first that $M$ is $2$-deep. The injectivity follows from \ref{chap11deep}, since $M$ is also $1$-deep. For the surjectivity, let $\psi$ be an element of $\mathrm{Hom}_A(I,M)$, for some admissible ideal $I = (f_1,\dots,f_r)$. Then $(\psi(f_i) f_i^{-1})_{1 \leq i \leq r}$ is a $1$-cycle in the \v{C}ech complex $C^{\bullet}(M,f_{\bullet})$. Since $M$ is $2$-deep, the module $H_I^1(M)$ vanishes by \ref{chap10.5bis}. By \ref{chap10.1}, the $1$-cycle $(\psi(f_i) f_i^{-1})_{1 \leq i \leq r}$ must be a $1$-boundary, so that $\psi(f_i) = f_i m$ in the localization of $M$ by $f_i$, for some $m$ in $M$. Thus $f_i^{N_i} \psi(f_i) = f_i^{N_i + 1} m $ in $M$ for some integers $N_i$. For $N = 1+ \sum_i N_i$ and $x$ in $I^{N}$ one has $\psi(x) = x m$, so that $x \psi(y)= y \psi(x) =xym$ for any $y$ in $I$. Since $M$ has no nonzero $I^N$-torsion this implies $\psi(y) = ym$, so that $\psi$ is the image of $m$.

Conversely, assume that the homomorphim
$$
M \rightarrow \mathrm{Hom}_A(I,M)
$$
is bijective for any admissible ideal $I$. The injectivity yields that $M$ is $1$-deep by \ref{chap11deep}. Let  $I = (f_1,\dots,f_r)$ be an admissible ideal. By \ref{chap10.1} we have to show that any $1$-cycle $(m_i f_i^{-N_i})_{1 \leq i \leq r}$ in the \v{C}ech complex $C^{\bullet}(M,f_{\bullet})$ is a $1$-boundary. By increasing the integers $N_i$ if necessary, one can assume that $f_i^{N_i} m_j = f_j^{N_j} m_i$ in $M$. For $N = \sum_i N_i$ and for any $r$-uple $\alpha = (\alpha_1,\dots,\alpha_r)$ of nonnegative integers summing to $N$, the element $m_{\alpha} = f^{\alpha - N_i e_i} m_i$ does not depend on the choice of an index $i$ such that $\alpha_i \geq N_i$. Here $f^{\alpha} = f_1^{\alpha_1}\dots f_r^{\alpha_r}$ and $e_i$ is the $i$-th basis vector.
Since $f^{\alpha} m_{\beta} = f^{\beta} m_{\alpha}$, one must have $m_{\alpha} = f^{\alpha} m$ for some $m$ by \ref{chap1interest} and by the hypothesis applied to $I^N$. Thus our initial $1$-cycle $(m_i f_i^{-N_i})_{1 \leq i \leq r}$ is the boundary of $m$.

\begin{cor}\label{chap12deepbis} Let $\Phi_0$ be a set of finitely generated ideals of $A$ generating the family of constructible supports of $A$ as in \ref{chap1exempleadique}$(i)$. An $A$-module $M$ is $2$-deep if and only if for any $I$ in $\Phi_0$ the homomorphism
$$
M \rightarrow \mathrm{Hom}_A(I,M)
$$
which sends an element $m$ to $(x \mapsto xm)$ is bijective.
\end{cor}

Assume indeed that the homomorphism
$$
M \rightarrow \mathrm{Hom}_A(I,M)
$$
is bijective for any $I$ in $\Phi_0$. It follows that for any $I$ in $\Phi_0$, $M$ has no nonzero $I$-torsion. However, if $I$ and $J$ are admissible ideal and if $M$ has no nonzero $I$-torsion then 
$$
\mathrm{Hom}_A(I,\mathrm{Hom}_A(J,M)) \cong \mathrm{Hom}_A(IJ,M),
$$
since the product homomorphism $I \otimes_A J \rightarrow IJ$ has $I$-torsion kernel. From this one deduces that the homomorphism
$$
M \rightarrow \mathrm{Hom}_A(I,M)
$$
is bijective whenever $I$ is a product of element of $\Phi_0$. The proof of Corollary \ref{chap12deep} then shows that $H^1_{I}(M)$ vanishes whenever $I$ is a product of elements of $\Phi_0$. One concludes with Corollary \ref{chap1comp}.

\subsection{\label{chap11.6}} In (\cite{EGA42}, $5.9$, $5.10$) the notions of ``\textit{modules purs}'' and ``\textit{modules clos}'' are introduced in a noetherian setting. We define here the corresponding purification and closure functors in a more general setting.

\begin{defi}\label{chap1defipurclos} Let $A$ be a $\Phi$-ring and let $M$ be an $A$-module. The \textbf{purification} of $M$ is defined as the quotient
$$
M^{\pur} = M / H_{\Phi_A}^0(M).
$$
If $M$ is $1$-deep (cf. \ref{chap10.5}), the \textbf{closure} of $M$ is defined as the colimit
$$
M^{\vartriangleleft} = \mathrm{colim} \ \mathrm{Hom}_{A}(I,M),
$$
where the colimit runs over all the admissible ideals of $A$. In general, we set $M^{\vartriangleleft} = (M^{\pur})^{\vartriangleleft}$ (cf. \ref{chap1purclos} below).
\end{defi}

\begin{exemple} Let $A$ be a valuation ring, and let us endow $A$ with the family of supports consisting of all the nonzero principal ideals of $A$. Then, for any $A$-module $M$, the purification $M^{\pur}$ is the largest torsion-free quotient of $M$, while we have
$$
M^{\vartriangleleft} = S^{-1} M,
$$
where $S = A \setminus \{ 0\}$.
\end{exemple}

\begin{prop}\label{chap1purclos} Let $A$ be a $\Phi$-ring and let $M$ be an $A$-module. Then $M^{\pur}$ is $1$-deep and $M^{\vartriangleleft}$ is $2$-deep.
\end{prop}

The $A$-module $M^{\pur}$ is $1$-deep by the characterization \ref{chap11deep}. In order to prove that $M^{\vartriangleleft}$ is $2$-deep, we can assume that $M$ is $1$-deep. Let $J$ be an admissible ideal of $A$. By \ref{chap1interest} the functor $\mathrm{Hom}_{A}(J,-)$ coincides on $1$-deep $A$-modules with $\mathrm{Hom}_{A}(G,-)$ for some finitely presented $A$-module $G$. In particular it commutes with filtered colimits of $1$-deep $A$-modules. Thus
\begin{align*}
\mathrm{Hom}_{A}(J,M^{\vartriangleleft}) &\cong \mathrm{colim} \ \mathrm{Hom}_{A}(J,\mathrm{Hom}_{A}(I,M)) \\
&\cong \mathrm{colim} \ \mathrm{Hom}_{A}(I \otimes_A J,M) \\
&\cong \mathrm{colim} \ \mathrm{Hom}_{A}(I J,M)\\
& \cong M^{\vartriangleleft}.
\end{align*}
The third identification above follows from the fact that the product map $I \otimes_A J \rightarrow IJ$ has $I$-torsion kernel (and $J$-torsion kernel as well). We conclude by Proposition \ref{chap12deep} that $M^{\vartriangleleft}$ is $2$-deep.

\begin{defi}\label{chap1defideep} The category $\Phi \mathrm{Rings}^{\geq d}$ is the full subcategory of $\Phi\mathrm{Rings}$ whose objects are the $\Phi$-rings which are $d$-deep as modules over themselves (cf. \ref{chap10.5}).
\end{defi}

There is a fully faithful functor from the category of rings to the category of $\Phi$-rings, which sends a ring $A$ to the $\Phi$-ring whose underlying ring is $A$ and whose only admissible ideal is the unit ideal. The image of this fully faithful functor is the intersection over all integers $d$ of the subcategories $\Phi\mathrm{Rings}^{\geq d}$ of $\Phi\mathrm{Rings}$. Indeed if for some ideal $I = (f_1,\dots,f_r)$ the \v{C}ech complex $C^{\bullet}(A,f_{\bullet})$ of Proposition \ref{chap10.1} is exact then it is an exact complex of flat $A$-modules so that $C^{\bullet}(A,f_{\bullet}) \otimes_A A / I = A / I [0]$ is exact as well, and thus $I = A$.

\begin{prop} \label{chap11.9} The \textbf{purification functor} 
$$
A \mapsto A^{\pur} = A / H_{\Phi_A}^0(A)
$$
from $\Phi \mathrm{Rings}$ to $\Phi \mathrm{Rings}^{\geq 1}$ yields a left adjoint to the canonical inclusion functor. Here $A / H_{\Phi_A}^0(A)$ is endowed with the family of constructible supports generated by that of $A$.
\end{prop}

This follows from the characterization of $1$-deep modules in \ref{chap11deep}.

\begin{prop} \label{chap11.9bis} The \textbf{closure functor}
$$
A \mapsto  A^{\vartriangleleft} = \mathrm{colim} \ \mathrm{Hom}_{A}(I,A),
$$
from $\Phi \mathrm{Rings}^{\geq 1}$ to $\Phi \mathrm{Rings}^{\geq 2}$ yields a left adjoint to the canonical inclusion functor. Here $A^{\vartriangleleft}$ is endowed with the family of constructible supports generated by that of $A$.
\end{prop}

Let $A$ be an object of $\Phi \mathrm{Rings}^{\geq 1}$. We notice as above that for any pair $I,J$ of admissible ideals, the product map $I \otimes_A J \rightarrow IJ$ has $I$-torsion kernel. This enables us to define an $A$-bilinear morphism
$$
\mathrm{Hom}_{A}(I,A) \times \mathrm{Hom}_{A}(J,A) \rightarrow \mathrm{Hom}_{A}(I \otimes_A J,A) \cong \mathrm{Hom}_{A}(IJ,A),
$$
which yields the $A$-algebra structure on $A^{\vartriangleleft}$. By Proposition \ref{chap1purclos}, the $A$-module $A^{\triangleleft}$ is $2$-deep. By Corollary \ref{chap1transitive}, we deduce that $A^{\triangleleft}$, when endowed with the family of constructible supports generated by that of $A$, is also $2$-deep as a module over itself. Using \ref{chap12deep} and \ref{chap1interest} again, one checks that any $\Phi$-ring morphism from $A$ to a $2$-deep $\Phi$-ring $B$ factors as an $A$-module homomorphism through $\mathrm{Hom}_{A}(I,A)$ for any admissible ideal $I$, hence through $A^{\vartriangleleft}$. The resulting homomorphism $A^{\vartriangleleft} \rightarrow B$ of $A$-modules must be a $\Phi$-ring morphism, and this shows that the morphism $A \rightarrow A^{\vartriangleleft}$ is universal among $\Phi$-ring morphisms from $A$ to a $2$-deep $\Phi$-ring.

%
%Let $J$ be an admissible ideal. By \ref{chap1interest} the functor $\mathrm{Hom}_{A}(J,-)$ coincides on $1$-deep $A$-modules with $\mathrm{Hom}_{A}(G,-)$ for some finitely presented $A$-module $G$. In particular it commutes with filtered colimits of $1$-deep $A$-modules. Thus
%\begin{align*}
%\mathrm{Hom}_{A}(J,A^{\vartriangleleft}) &\xleftarrow{\sim} \mathrm{colim} \ \mathrm{Hom}_{A}(J,\mathrm{Hom}_{A}(I,A)) \\
%&\xleftarrow{\sim} \mathrm{colim} \ \mathrm{Hom}_{A}(IJ,A) = A^{\vartriangleleft}.
%\end{align*}
%Since the left hand side is equal to $\mathrm{Hom}_{A^{\vartriangleleft}}(JA^{\vartriangleleft},A^{\vartriangleleft})$ by \ref{chap1interest} again, we conclude by Corollary \ref{chap12deepbis} that $A^{\vartriangleleft}$ is $2$-deep when endowed with the family of constructible supports generated by that of $A$. Using \ref{chap1interest} again, one checks that any $\Phi$-ring morphism from $A$ to a $2$-deep $\Phi$-ring $B$ factors as an $A$-module homomorphism through $\mathrm{Hom}_{A}(I,A)$ for any admissible ideal $I$, hence through $A^{\vartriangleleft}$. The resulting homomorphism $A^{\vartriangleleft} \rightarrow B$ of $A$-modules must be a $\Phi$-ring morphism, and this shows that the morphism $A \rightarrow A^{\vartriangleleft}$ is universal among $\Phi$-ring morphisms from $A$ to a $2$-deep $\Phi$-ring.

The composition of the functors of propositions \ref{chap11.9} and \ref{chap11.9bis} will still be denoted by $A \rightarrow A^{\vartriangleleft} $.

\begin{defi} \label{chap1vringdef} The category $\Phi \mathrm{Rings}^+$ of \textbf{normal $\Phi$-rings} is the full subcategory of $\Phi \mathrm{Rings}$ whose objects are the $\Phi$-rings $A$ such that the adjunction morphism
$$
A \rightarrow A^{\vartriangleleft}
$$
is injective, so that $A$ can be considered as a subring of $A^{\vartriangleleft}$, and such that $A$ is integrally closed in $A^{\vartriangleleft}$: any element of $A^{\vartriangleleft}$ which is integral over $A$ is required to lie in $A$.
\end{defi}

In particular, any object of $\Phi \mathrm{Rings}^{\geq 2}$ is an object of $\Phi \mathrm{Rings}^+$, and any object of $\Phi \mathrm{Rings}^+$ is an object in $\Phi \mathrm{Rings}^{\geq 1}$. 

\begin{prop} \label{chap1110} The \textbf{normalization functor} which sends a $\Phi$-ring $A$ to the integral closure $A^+$ of $A$ in its closure $A^{\vartriangleleft}$ yields a left adjoint to the canonical inclusion functor from $\Phi \mathrm{Rings}^+$ to $\Phi \mathrm{Rings}$. Here $A^+$ is endowed with the family of constructible supports generated by that of $A$.
\end{prop}

Let $f: A \rightarrow B$ be a morphism in $\Phi \mathrm{Rings}$ such that $B$ belongs to $\Phi \mathrm{Rings}^+$. This induces a morphism $A^{\vartriangleleft} \rightarrow B^{\vartriangleleft}$ sending the image of $A$ into $B$. Since $B \subseteq B^{\vartriangleleft}$ is an integrally closed subring, the morphism $A \rightarrow B$ factors uniquely through $A \rightarrow A^+$, the integral closure of the image of $A$ in $A^{\vartriangleleft}$. Thus the morphism $A \rightarrow A^+$ is universal among $\Phi$-ring morphisms from $A$ to a normal $\Phi$-ring.

\subsection{\label{chap1111}} We now turn to the definition and study of the local objects in our category of $\Phi$-rings.

\begin{defi}\label{chap1defphilocal} A $\Phi$-ring $A$ is \textbf{$\Phi$-local} if it is a local ring and if all of its admissible ideals are invertible.
\end{defi}

In (\cite{Abbes11} 1.9.1), these objects are called ``prevaluative'' rings, in the particular case where the admissible ideals are the open finitely generated ideals for some adic topology. 

\begin{exemple} Let $A$ be a valuation ring, and let us endow $A$ with the family of supports consisting of all the nonzero principal ideals of $A$. Then $A$ is a $\Phi$-local $\Phi$-ring.
\end{exemple}

\begin{rema}\label{chap1LF1} Let $A$ be a $\Phi$-local $\Phi$-ring. The condition that any admissible ideal of $A$ is invertible implies that whenever an admissible ideal $I$ of $A$ is generated by a family $(x_\lambda)_{\lambda \in \Lambda}$ of elements of $A$, there exists $\lambda$ in $\Lambda$ such that $x_\lambda$ is a nonzero divisor generating $I$. Indeed, the ideal $I$ is an invertible ideal in a local ring, hence is generated by some nonzero divisor $x$ of $A$. For each $\lambda$ in $\Lambda$, we then have $x_{\lambda} = a_{\lambda} x$ for a unique element $a_{\lambda}$ of $A$. Moreover, we can write $x$ as $\sum_{\lambda \in \Lambda} b_{\lambda} x_{\lambda}$ for some family $(b_\lambda)_{\lambda \in \Lambda}$ of elements of $A$, only finitely many of which being non zero. Since $x$ is a non zero divisor, we must have $\sum_{\lambda \in \Lambda} a_{\lambda} b_{\lambda} = 1$. Since $A$ is local, there exists $\lambda$ in $\Lambda$ such that $a_{\lambda}$ is a unit, in which case $x_\lambda$ is a nonzero divisor generating $I$.
\end{rema}

A structure theorem similar to (\cite{Abbes11} 1.9.4) holds in the general case:

\begin{prop}\label{chap12.1} Let $A$ be a $\Phi$-local $\Phi$-ring with maximal ideal $\mathfrak{n}$.
\begin{itemize}
\item[$(i)$] The canonical homomorphism $A \rightarrow A^{\triangleleft}$ is injective.
\item[$(ii)$] The ring $A^{\vartriangleleft}$ is a local ring, and its maximal ideal $\mathfrak{m}$ is contained in $A$. More precisely, $A^{\vartriangleleft}$ is the localization of $A$ at $\mathfrak{m}$.
\item[$(iii)$] An ideal of $A$ is admissible if and only if it is generated by an element of $A \setminus \mathfrak{m}$.
\item[$(iv)$] The ring $A/\mathfrak{m}$ is a valuation subring of the field $A^{\vartriangleleft}/\mathfrak{m}$.
\item[$(v)$] We have
$$
\mathfrak{m} = \bigcap_{s \in A \setminus \mathfrak{m}} s \mathfrak{n}.
$$
\item[$(vi)$] $A$ is a normal $\Phi$-ring.
\item[$(vii)$] Let $I$ be a finitely generated ideal of $A$ and let $g$ be an element of $A$ such that $(I,g)$ is admissible. Then $I (A/\mathfrak{m}) \subseteq g(A/\mathfrak{m})$ if and only if $I \subseteq gA$.
\end{itemize}

\end{prop}

Let $S$ be the subset of elements of $A$ which generate an admissible ideal. Then $S$ is a multiplicative subset of the set of nonzero divisors of $A$ and any admissible ideal of $A$ is equal to $s A$ for some element $s$ of $S$. In particular we have
$$
A^{\triangleleft} \cong \underset{s \in S}{\mathrm{colim}} \ \mathrm{Hom}_A(sA,A) \xrightarrow{\sim}  \underset{s \in S}{\mathrm{colim}} \ s^{-1}A \xrightarrow{\sim} S^{-1} A.
$$
This proves in particular $(i)$. Let us define
$$
\mathfrak{m} = \bigcap_{s \in S} s \mathfrak{n}.
$$
Then $\mathfrak{m}$ is a proper ideal of $A^{\vartriangleleft}$ contained in $A$. We show that $S = A \setminus \mathfrak{m} $.

First, let $a$ be an element of $A \setminus \mathfrak{m}$. Then there exists an element $s$ of $S$ such that $a$ does not belong to $s \mathfrak{n}$. The finitely generated ideal $(s,a)$ is admissible, as it contains the admissible ideal $sA$. Thus $(s,a)$ is invertible, since $A$ is $\Phi$-local, and is in particular generated by either $a$ or $s$, cf. \ref{chap1LF1}.
If $(s,a)$ is generated by $s$, then $a = bs$ for some $b$ in $A$. But then $b$ does not belong to $\mathfrak{n}$, so it is a unit of $A$, and hence $s = b^{-1}a$ belongs to $aA$.
If $(s,a)$ is generated by $a$, it is also the case that $s$ belongs to $aA$. Thus $aA$ contains the admissible ideal $sA$ and is therefore admissible itself, so that $a$ belongs to $S$. 
Conversely, if $a$ is an element of $S$, then $a$ does not belong to $a \mathfrak{n}$, and thus $a$ is in $A \setminus \mathfrak{m}$. 

This shows $(iii)$ and that $A^{\vartriangleleft}$ is the localization of $A$ at $\mathfrak{m}$, hence $(ii)$ and $(v)$. Moreover, for any pair $a,s$ of elements of $A \setminus \mathfrak{m}$, the ideal $(s,a)$ of $A$ is admissible, hence it is generated by either $a$ or $s$, cf. \ref{chap1LF1}. Thus, either $a$ belongs to $sA$ or $s$ belongs to $aA$. This shows that $A/\mathfrak{m}$ is a valuation subring of the field $A^{\vartriangleleft}/\mathfrak{m}$, hence $(iv)$.

In particular, $A/\mathfrak{m}$ is integrally closed in $A^{\vartriangleleft}/\mathfrak{m}$, so that $A$ is integrally closed in $A^{\vartriangleleft}$. Thus $A$ is a normal $\Phi$-ring, hence $(vi)$.

Let $I$ be a finitely generated ideal of $A$ and let $g$ be an element of $A$ such that $(I,g)$ is admissible, hence equal to $sA$ for some element $s$ of $S$. If we have the inclusion $I (A/\mathfrak{m}) \subseteq g(A/\mathfrak{m})$ then $s (A/\mathfrak{m}) \subseteq g(A/\mathfrak{m})$ also holds, so that $s = ag + b$ with $a \in A$ and $b \in \mathfrak{m}$. Since $\mathfrak{m}\subseteq s \mathfrak{n}$ we have $b = c s$ for some $c \in \mathfrak{n}$. Thus we have $(1-c) s = ag$, and consequently $s = (1-c)^{-1} a g$, so that $(I,g) = sA$ is contained in $gA$. In particular, we have $I \subseteq gA$, hence $(vii)$.

\begin{rema}\label{chap1remarkphilocal} The datum of a $\Phi$-local $\Phi$-ring is the same as the datum of a local ring together with a valuation subring of its residue field. The inverse construction is simply the following: given a local ring $B$ and a valuation subring $R$ of its residue field, let $A$ be the inverse image of $R$ in $B$, and declare an ideal of $A$ to be admissible if it is generated by an element which is invertible in $B$. Then $A$ is a $\Phi$-ring and $A^{\vartriangleleft} = B$. Moreover, the valuation subring of the residue field of $A^{\vartriangleleft}$ produced by Proposition \ref{chap12.1} is precisely $R$.
\end{rema}

\begin{constr}\label{chap1Constr}Let $A$ be a local $\Phi$-ring and let $f : A \rightarrow S$ be a local homomorphism from $A$ to a valuation ring $S$ such that $f(I)S \neq 0$ for any admissible ideal $I$ of $A$. Let us endow $S$ with the family of constructible supports generated by that of $A$. Then $S$ is a $\Phi$-local $\Phi$-ring. Let $\mathfrak{m}$ be the maximal ideal of $S^{\triangleleft}$ (cf. \ref{chap12.1}) and let $\p$ be its inverse image in $A^{\triangleleft}$. The homomorphism $A^{\triangleleft} \rightarrow S^{\triangleleft}/\mathfrak{m}$ induced by $f$ factors through an injective homomorphism from $\ka(\p) = \mathrm{Frac}(A^{\triangleleft}/ \p)$ to $S^{\triangleleft}/\mathfrak{m}$. Let $R \subseteq \ka(\p)$ be the inverse image of the valuation ring $S/\mathfrak{m}$ under this homomorphism. The pair $(\p,R)$ has the following properties:
\begin{itemize}
\item[(a)] The ring $R$ is a valuation ring, and the homomorphism $A \rightarrow \ka(\p)$ factors through a local homomorphism $A \rightarrow R$.
\item[(b)] For any $a \in \ka(\p)$ there is an admissible ideal $I$ of $A$ such that $Ia \subseteq R$.
\item[(c)] The prime ideal $\p$ contains no admissible ideal of $A^{\triangleleft}$.
\item[(d)] Let $I$ be a finitely generated ideal of $A$ and let $g$ be an element of $A$ such that $(I,g)$ is admissible. Then $IR \subseteq gR$ if and only if $IS \subseteq gS$.
\end{itemize}

If $a$ is a non zero element of $R$, then $a$ or $a^{-1}$ belongs to the valuation ring $S/\mathfrak{m}$, hence to $R$. Thus $R$ is a valuation ring. Since its maximal ideal is the inverse image of that of $S/\mathfrak{m}$, the local homomorphism $A \rightarrow S/\mathfrak{m}$ factors through a local homomorphism $A \rightarrow R$.

For any $a \in \ka(\p)$, there is an admissible ideal $I$ of $A$ such that $f(I a)(S/\mathfrak{m}) \subseteq S/\mathfrak{m}$. For such an ideal $I$, we have $I a \subseteq R$.

The maximal ideal $\mathfrak{m}$ of $S^{\triangleleft}$ is the intersection of all $f(I) \mathfrak{n}$ where $I$ runs over the admissible ideals of $A$, and where $\mathfrak{n}$ is the maximal ideal of $S$. For any admissible ideal $I$ of $A$, the invertible ideal $f(I)S$ is not contained in $f(I) \mathfrak{n}$, and is consequently not contained in $\mathfrak{m}$. In particular, the ideal $I$ is not contained in $\p$.

 Let $I$ be a finitely generated ideal of $A$ and let $g$ be an element of $A$ such that $(I,g)$ is admissible. Then $IR \subseteq gR$ if and only if $I(S/\mathfrak{m}) \subseteq g (S/\mathfrak{m})$, if and only if $IS \subseteq gS$ by \ref{chap12.1} (vii).
\end{constr}

\section{Sheaves of $\Phi$-rings \label{chap1phiringssheaves}}

Let $\mathcal{U}$ be a universe (\cite{SGA4}, I.0). For any category $C$, we denote by $\mathcal{U}\text{-}C$ the full subcategory of $C$ whose objects are the objects of $C$ which belong to $\mathcal{U}$.

\subsection{\label{chap10.2.0}} Let $C$ be a $\mathcal{U}$-site (\cite{SGA4}, II.3.0.2). A \textbf{presheaf of $\Phi$-rings} on $C$ (with respect to $\mathcal{U}$) is a contravariant functor from $C$ to $\mathcal{U}\text{-}\Phi \mathrm{Rings}$. A presheaf of $\Phi$-rings $A$ (with respect to $\mathcal{U}$) is a \textbf{sheaf of $\Phi$-rings} if for any $\Phi$-ring $B$ in $\mathcal{U}$, the functor
$$
U \in C \mapsto \mathrm{Hom}_{\Phi \mathrm{Rings}}(B,A(U))
$$
is a sheaf of sets on $C$ (with respect to $\mathcal{U}$). We also define a morphism of sheaves of $\Phi$-rings to be a morphism of presheaves of $\Phi$-rings between sheaves of $\Phi$-rings. We henceforth omit references to $\mathcal{U}$ when no confusion arise from the lack thereof.

This definition is a specialization of (\cite{SGA4}, II.6.1).  Since the category $\Phi \mathrm{Rings}$ is complete \ref{chap1compcocomp}, a presheaf of $\Phi$-rings on $C$ is a sheaf of $\Phi$-rings if and only if for any covering sieve $S$ of an object $U$ of $C$, the natural morphism
\begin{align}\label{chap1eq1}
A(U) \rightarrow \lim_{C_{/S}} A
\end{align}
is an isomorphism of $\Phi$-rings (\cite{SGA4}, II.6.2).

\begin{prop}\label{chap1sheafcond} Let $A$ be a presheaf of $\Phi$-rings on $ C$. Then $A$ is a sheaf of $\Phi$-rings if and only if its underlying presheaf of sets is a sheaf and for each integer $r$ the contravariant functor 
$$
\mathrm{Adm}_r(A) : U \mapsto \{ (f_i)_{1 \leq i \leq r} \in A(U)^r \ | \ (f_i)_{1 \leq i \leq r} \text{ generates an admissible ideal of } A(U) \}
$$
is a sheaf of sets on $C$.
\end{prop}

Indeed if $A$ is a sheaf of $\Phi$-rings then its underlying presheaf of rings is a sheaf since the forgetful functor from $ \Phi \mathrm{Rings}$ to the category of rings has a left adjoint by \ref{chap1commute} and thus commutes with the limit in \ref{chap1eq1}. For each integer $r$, taking the $\Phi$-ring $B$ in the definition of a sheaf of $\Phi$-rings to be the polynomial ring $\mathbb{Z}[X_1,\dots,X_r]$ endowed with the family of constructible supports generated by the ideal $(X_1,\dots,X_r)$ yields that the presheaf $\mathrm{Adm}_r(A)$ is a sheaf of sets.

Conversely, let us assume that the presheaf of rings underlying $A$ is a sheaf and that $\mathrm{Adm}_r(A)$ is a sheaf for each integer $r$. For each covering sieve $S$ of an object $U$ of $C$, the natural homomorphism
\begin{align}\label{chap1eq2}
 A(U) \rightarrow \lim_{C_{/ S}} A
\end{align}
is bijective since $A$ is a sheaf of rings. Moreover an element $(f_i)_{1 \leq i \leq r}$ of  $A(U)^r$ generates an admissible ideal if and only if it generates an admissible ideal of $A(V)^r$ for each $V \rightarrow U$ in $S$, since $\mathrm{Adm}_r(A)$ is a subsheaf of $A^r$. According to the description of limits in $ \Phi \mathrm{Rings}$ given in the proof of \ref{chap1compcocomp} the latter condition holds if and only if $(f_i)$ generates an admissible ideal of $\lim_{C_{/ S}} A$. Thus \ref{chap1eq2} is an isomorphims of $\Phi$-rings and $A$ is a sheaf of $\Phi$-rings.

\begin{cor}\label{chap1sheafcond2} Let $A$ be a presheaf of $\Phi$-rings on $C$ and let $\widetilde{A}$ be the sheaf of rings associated to its underlying presheaf of rings. Then there is a unique sheaf of $\Phi$-rings with underlying sheaf of rings $\widetilde{A}$, such that for each $r$ the sheaf $\mathrm{Adm}_r(\widetilde{A})$ (cf. \ref{chap1sheafcond}) is the sheaf associated to $\mathrm{Adm}_r(A)$. Moreover the morphism $A \rightarrow \widetilde{A}$ is initial among the morphisms of presheaves of $\Phi$-rings from $A$ to a sheaf of $\Phi$-rings.
\end{cor}

Let $\mathrm{Adm}_r(\widetilde{A}) \subseteq \widetilde{A}^r$ be the sheaf associated to $\mathrm{Adm}_r(A)$.
For any object $U$ of $C$, let us declare an ideal of $\widetilde{A}(U)$ admissible if it is generated by elements $f_1,\dots,f_r$ such that  $(f_1,\dots,f_r)$ belongs to $\mathrm{Adm}_r(\widetilde{A})(U)$.

Let $r$ and $s$ be integers and let us consider the morphism of presheaves of sets given by
\begin{align*}
\alpha_{r,s} : A^{rs} \times A^s &\rightarrow A^r \\
((a_{ij})_{1 \leq i \leq r, 1 \leq j \leq s}, (g_j)_{1 \leq j \leq s}) &\mapsto \left( \sum_{j} a_{ij} g_j \right)_{1 \leq i \leq r}.
\end{align*}
The fiber product $\alpha_{r,s}^{-1}(\mathrm{Adm}_r(A)) = (A^{rs} \times A^s) \times_{A^r} \mathrm{Adm}_r(A)$ is a subobject of $A^{rs} \times A^s$ contained in $A^{rs} \times \mathrm{Adm}_s(A)$ since $A$ is a presheaf of $\Phi$-rings (cf. \ref{chap111}$(2)$). Since the ``associated sheaf'' functor (\cite{SGA4}, II.3.4) commutes with finite limits (\cite{SGA4}, II.4.1), we obtain that the morphism $\widetilde{\alpha}_{r,s}$ between the associated sheaves of sets satisfies $\widetilde{\alpha}_{r,s}^{-1}(\mathrm{Adm}_r(\widetilde{A})) \subseteq \widetilde{A}^{rs} \times \mathrm{Adm}_s(\widetilde{A})$. If $U$ is an object of $C$ then this shows that the axiom \ref{chap111}$(2)$ is satisfied for $\widetilde{A}(U)$ and that an element $(f_1,\dots,f_r)$ of $\widetilde{A}^r(U)$ generates an admissible ideal if and only if it belongs to $\mathrm{Adm}_r(\widetilde{A})(U)$. 

By considering
\begin{align*}
\beta_{r,s} : \mathrm{Adm}_r(A) \times \mathrm{Adm}_s(A) &\rightarrow \mathrm{Adm}_{rs}(A) \\
((f_i)_{1 \leq i \leq r}, (g_j)_{1 \leq j \leq s}) &\mapsto \left( f_i g_j \right)_{1 \leq i \leq r,1 \leq j \leq s}
\end{align*}
instead of $\alpha_{r,s}$ and by taking the induced morphism between the associated sheaves of sets, one obtains the axiom \ref{chap111}$(1)$ as well. Thus we have endowed the sheaf of rings $\widetilde{A}$ with a structure of presheaf of $\Phi$-rings. Since a local section $(f_1,\dots,f_r)$ of $\widetilde{A}^r$ generates an admissible ideal if and only if it belongs to $\mathrm{Adm}_r(\widetilde{A})$, the notation $\mathrm{Adm}_r(\widetilde{A})$ is consistent with \ref{chap1sheafcond}. For each $r$, $\mathrm{Adm}_r(\widetilde{A})$ is a sheaf, hence $\widetilde{A}$ is a sheaf of $\Phi$-rings by Proposition \ref{chap1sheafcond}.

Let $A \rightarrow B$ be a morphism of presheaves of $\Phi$-rings from $A$ to a sheaf of $\Phi$-rings. The underlying morphism of presheaves of rings uniquely factors through the natural morphism $A \rightarrow \widetilde{A}$. Since for each integer the morphism $A^r \rightarrow B^r$ maps $\mathrm{Adm}_r(A)$ into $\mathrm{Adm}_r(B)$, the induced morphism $\widetilde{A}^r \rightarrow B^r$ maps $\mathrm{Adm}_r(\widetilde{A})$ into $\mathrm{Adm}_r(B)$ as well. Thus the morphism $\widetilde{A} \rightarrow B$ is a morphism of sheaves of $\Phi$-rings.

\begin{cor}\label{chap1sheafcond3} Let $A$ be a sheaf of $\Phi$-rings on $C$, and let $U$ be a sheaf of sets on $C$. Then the ring $\mathrm{Hom}(U,A)$ is endowed with a family of constructible supports by declaring admissible any ideal of $\mathrm{Hom}(U,A)$ which is generated by elements $f_1,\dots,f_r$ such that $(f_1,\dots,f_r)$ belongs to the subset $\mathrm{Hom}(U,\mathrm{Adm}_r(A))$ of $\mathrm{Hom}(U,A^r)=\mathrm{Hom}(U,A)^r$.
\end{cor}

This follows from the proof of \ref{chap1sheafcond2} by applying the limit preserving functor $\mathrm{Hom}(U,\cdot)$ to $\alpha_{r,s}$ and $\beta_{r,s}$. 

\subsection{\label{chap10.2.0.1}} Let $X$ be a topos. We define a \textbf{$\Phi$-ring of $X$} to be a ring $A$ in $X$ together with a subobject $\mathrm{Adm}_r(A)$ of $A^r$ for each positive integer $r$, such that:
\begin{itemize}
\item[$(1)$] For all $U$ in $X$, the element $1 \in A(U)$ belongs to $\mathrm{Adm}_1(A)(U)$. For all positive integers $r,s$, for all $(f_i)_{1 \leq i \leq r}$ and $(g_j)_{1 \leq j \leq s}$ in $\mathrm{Adm}_r(A)(U)$ and $\mathrm{Adm}_s(A)(U)$ respectively, the element $\left( f_i g_j \right)_{1 \leq i \leq r,1 \leq j \leq s}$ of $A(U)^{rs}$ belongs to $\mathrm{Adm}_{rs}(A)(U)$. 
\item[$(2)$] For all $U$ in $X$, for all positive integers $r,s,$ and for all $(f_i)_{1 \leq i \leq r}$ in $\mathrm{Adm}_r(A)(U)$, any element $(g_j)_{1 \leq j \leq s}$ of $A(U)^s$ such that $\sum_i f_i A(U) \subseteq \sum_j g_j A(U)$ belongs to $\mathrm{Adm}_{s}(A)(U)$. 
\end{itemize}
We define a morphism of $\Phi$-rings of $X$ to be a ring homomorphism $f : A \rightarrow B$ such that for each positive integer $r$, the subobject $f(\mathrm{Adm}_r(A))$ of $B^r$ is a subobject of $\mathrm{Adm}_r(B)$. This defines a category $\Phi \mathrm{Rings}_X$. 

\begin{prop}\label{chap1sheafphicomp} Let $C$ be a site and let $\widetilde{C}$ be the topos of sheaves of sets on $C$. The functor which associates to a sheaf of $\Phi$-rings $A$ on $C$ the $\Phi$-ring $(A,(\mathrm{Adm}_r(A))_{r\geq 0})$ in $\widetilde{C}$ (cf. \ref{chap1sheafcond}) is an equivalence of categories.
\end{prop}

This follows from \ref{chap1sheafcond} and \ref{chap1sheafcond3}.

\subsection{\label{chap10.2.0.2}} Let $f = (f_*,f^{-1}) : X \rightarrow Y$ be a morphism of toposes\footnote{We choose here the notation $f^{-1}$ in order to avoid confusion with the pullback functor $f^*$ between categories of modules, which is associated to a morphism a ringed toposes.}. If $A$ is a $\Phi$-ring of $X$ then setting $\mathrm{Adm}_r(f_* A) = f_* (\mathrm{Adm}_r(A))$ gives $f_* A$ the structure of a $\Phi$-ring of $Y$. Similarly if $B$ is a $\Phi$-ring of $X$ then setting $\mathrm{Adm}_r(f^{-1} B) = f^{-1} \mathrm{Adm}_r(B)$ gives $f^{-1} B$ the structure of a $\Phi$-ring of $Y$, since the functor $f^{-1}$ commutes with finite limits.

In particular if $A$ is a $\Phi$-ring of $X$ and $x$ is a point of $X$, then \textbf{the stalk $A_x$ of $A$ at $x$} is defined to be $x^{-1} A$. It is a $\Phi$-ring of the topos of sets, or in other words a $\Phi$-ring. We have 
$$
A_x = \mathrm{colim}_{(U,u)} A(U),
$$
in the category of $\Phi$-rings, where the colimit is indexed by pairs $(U,u)$ where $U$ runs over the members of a small generating family of $X$ and $u$ is an element of $U_x = x^{-1} U$.

\begin{prop} The pair of functors $(f_*,f^{-1})$ between the categories $\Phi \mathrm{Rings}_X$ and $\Phi \mathrm{Rings}_Y$ is a pair of adjoint functors.
\end{prop}

This follows from the corresponding adjunction property of $(f_*,f^{-1})$ between the categories $X$ and $Y$.

\subsection{\label{chap10.2.0.2.1}} A $\Phi$-ringed topos is a pair $(X,\Ow_X)$ where $X$ is a topos and $\Ow_X$ is a $\Phi$-ring of $X$ (cf. \ref{chap10.2.0.1}). We define a morphism of $\Phi$-ringed topos from a $\Phi$-ringed topos $(X,\Ow_X)$ to an other one $(Y,\Ow_Y)$ to be a morphism of toposes $f = (f_*,f^{-1}) : X \rightarrow Y$ together with a morphism $f^{\sharp} : \Ow_Y \rightarrow f_* \Ow_X$ (cf. \ref{chap10.2.0.2}) of $\Phi$-rings in $Y$.

\subsection{\label{chap10.2.0.3}} Let us recall that a ring $A$ in a topos $X$ is \textbf{local} if the natural morphism
$$
A^{\times} \sqcup (1+A^{\times}) \rightarrow A
$$
is an epimorphism and if the limit of the diagram $e \xrightarrow[]{1} A \xleftarrow[]{0} e$ is the initial object of $X$\footnote{The last condition, which can be reworded as ``$0 \neq 1$'', is missing from the definition of a local ring in a topos which can be found in (\cite{SGA4}, IV.13.9)}.

\begin{prop}\label{chap1localringed} Let $(X,\Ow_X)$ be a ringed topos with enough points. Then $\Ow_X$ is local if and only if for each point $x$ of $X$ the stalk $\Ow_{X,x}$ of $\Ow_X$ at $x$ (cf. \ref{chap10.2.0.2}) is a local ring.
\end{prop}

This appears as an exercise in (\cite{SGA4}, IV.13.9). The proposition follows from the fact that a ring $A$ (in the punctual topos) is local if and only if for any element $a$ of $A$, either $a$ or $1+a$ is invertible, and if moreover $0 \neq 1$ in $A$. 

\begin{defi}\label{chap1localphiringed} A $\Phi$-ringed topos $(X,\Ow_X)$ (cf. \ref{chap10.2.0.2.1}) is locally $\Phi$-ringed if the underlying ring of $\Ow_X$ is local. A morphism of locally $\Phi$-ringed toposes is a morphism of $\Phi$-ringed toposes (cf. \ref{chap10.2.0.2.1}) which is also a morphism of locally ringed toposes (\cite{SGA4}, IV.13.9).
\end{defi}

\subsection{\label{chap10.2.0.3.5}} A $\Phi$-ring $A$ in a topos is $\Phi$-local if the following three conditions are satisfied:
\begin{itemize}
\item[$(1)$] The ring $A$ is local (cf. \ref{chap10.2.0.3}).
\item[$(2)$] For all $U$ in $X$, each element of $\mathrm{Adm}_1(A)(U)$ is not a zero divisor in $A(U)$.
\item[$(3)$] For each integer $r$ the morphism
\begin{align*}
\coprod_{j=1}^r \left( \mathrm{Adm}_1(A) \times A^{[1,r] \setminus \{j\}} \right) &\rightarrow \mathrm{Adm}_r(A) \\
(f, (a_{i})_{1 \leq i \leq r, i \neq j}) &\mapsto (a_i f)_{1 \leq i \leq r}  \text{ where } a_j = 1
\end{align*}
is epimorphic.
\end{itemize}

\begin{defi}\label{chap1defiphilocalringed} A $\Phi$-ringed topos $(X,\Ow_X)$ is $\Phi$-locally $\Phi$-ringed if $\Ow_X$ is $\Phi$-local. A morphism of $\Phi$-locally $\Phi$-ringed toposes is a morphism of locally $\Phi$-ringed toposes (cf. \ref{chap1localphiringed}) between $\Phi$-locally $\Phi$-ringed toposes.
\end{defi}

\begin{prop}\label{chap1philocalringed} Let $(X,\Ow_X)$ be a $\Phi$-ringed topos with enough points. Then $(X,\Ow_X)$ is $\Phi$-locally $\Phi$-ringed if and only if for each point $x$ of $X$ the stalk $\Ow_{X,x}$ of $\Ow_X$ at $x$ (cf. \ref{chap10.2.0.2}) is a $\Phi$-local $\Phi$-ring.
\end{prop}

This follows from Proposition \ref{chap1localringed} and from the fact that a $\Phi$-ring (in the punctual topos) is $\Phi$-local if and only if its underlying ring is local and if for any admissible ideal $I=(f_1,\dots,f_r)$ of $A$, there exists $i$ such that $f_i$ generates $I$ and is not a zero divisor.

\subsection{\label{chap10.2.0.4}} Let us recall that we denote by $\Phi \mathrm{Rings}_X$ the category of $\Phi$-rings of a topos $X$ (cf. \ref{chap10.2.0.1}).

\begin{defi} Let $(X,\Ow_X)$ be a $\Phi$-ringed topos. For $d \leq 2$, an $\Ow_X$-module $\M$ is said to be $d$-deep if for all $U$ in $X$, the $\Ow_X(U)$-module $\M(U)$ is $d$-deep (cf. \ref{chap10.5}).
\end{defi}

\begin{prop}\label{chap1deepstalk} Let $d \leq 2$ and let $(X,\Ow_X)$ be a ringed topos with enough points. Then an $\Ow_X$-module $\M$ is $d$-deep if and only if all of its stalks are $d$-deep.
\end{prop}

The case $d=0$ is tautological. By \ref{chap11deep} an $\Ow_X$-module $\M$ is $1$-deep if and only if for each integer $r$, the inverse image of the neutral suboject $0 \rightarrow \M^r$ of $\M^r$ by the morphism
\begin{align*}
\mathrm{Adm}_r(\Ow_X) \times \M &\rightarrow \M^r \\
((f_j)_{1 \leq j \leq r},m) &\mapsto (f_j m)_{1 \leq j \leq r}
\end{align*}
is a subobject of $\mathrm{Adm}_r(\Ow_X) \times 0$. This yields the proposition for $d=1$. For $d=2$, one rather uses \ref{chap1interest} and \ref{chap12deep}, which imply that a $\M$ is $2$-deep if and only if for each integer $r$ the natural morphism
\begin{align*}
\mathrm{Adm}_r(\Ow_X) \times \M &\rightarrow \mathrm{lim} \left( \mathrm{Adm}_r(\Ow_X) \times \M^r \rightrightarrows \M^{r^2} \right) \\
((f_j)_{1 \leq j \leq r},m) &\mapsto ((f_j)_{1 \leq j \leq r},(f_j m)_{1 \leq j \leq r})
\end{align*}
is an isomorphism, where the transition morphisms in the limit are given by $((f_j)_{ j },(m_j)_{j}) \mapsto (f_i m_j)_{i,j}$ and $((f_j)_{j},(m_j)_{j}) \mapsto (f_j m_i)_{i,j}$.

\begin{defi} Let $X$ be a topos. For $d \leq 2$ we define the category $\Phi \mathrm{Rings}_X^{\geq d}$ of $d$-deep $\Phi$-rings to be the full subcategory of $\Phi \mathrm{Rings}_X$ whose objects are the $\Phi$-rings $A$ which are $d$-deep as modules over themselves.
\end{defi}

\begin{prop}\label{chap1deepsheaf} Let $d \leq 2$ be an integer, and let $A$ be a presheaf of $\Phi$-rings on a site $C$ such that $A(U)$ is $d$-deep for any object $U$ of $C$. Then the sheaf of $\Phi$-rings associated to $A$ (cf. \ref{chap1sheafcond2}) is a $d$-deep $\Phi$-ring in the topos of sheaves of sets on $C$.
\end{prop}

This follows from the characterizations of $d$-deep modules in a $\Phi$-ringed topos given in the proof of the proposition \ref{chap1deepstalk}.

\begin{prop}\label{chap1sheafpurif} Let $X$ be a topos. The canonical inclusion functor from $\Phi \mathrm{Rings}_X^{\geq 1}$ to $\Phi \mathrm{Rings}_X$ admits a left adjoint. We call this left adjoint the \textbf{purification} functor and we denote it by $A \mapsto A^{\mathrm{pur}}$.
\end{prop}

Let $C$ be a site such that $X$ is (equivalent to) the topos of sheaves of sets on $C$ (e.g. the canonical site of $X$). By \ref{chap1sheafphicomp} the category $\Phi \mathrm{Rings}_X$ is isomorphic to the category of sheaves of $\Phi$-rings on $C$. A sheaf of $\Phi$-rings $A$ on $C$ corresponds to a $1$-deep $\Phi$-ring object in $X$ if and only if for each object $U$ of $C$ the $\Phi$-ring $A(U)$ is $1$-deep. Thus the functor which associated to a sheaf of $\Phi$-rings $A$ on $C$ the sheaf associates to the presheaf of $\Phi$-rings $U \mapsto A(U)^{\mathrm{pur}}$ (cf. \ref{chap1sheafcond2} and \ref{chap1deepsheaf}) defines a left adjoint to the canonical inclusion functor.

\begin{prop}\label{chap1sheafclos} Let $X$ be a topos. The canonical inclusion functor from $\Phi \mathrm{Rings}_X^{\geq 2}$ to $\Phi \mathrm{Rings}_X$ admits a left adjoint. We call it the \textbf{closure} functor and we denote it by $A \mapsto A^{\triangleleft}$.
\end{prop}

This is proved similarly as \ref{chap1sheafpurif} using the closure functor on $\Phi$-rings instead of the purification functor.

\begin{prop}\label{chap1sheafstalk} Let $X$ be a topos, let $x$ be a point of $X$ and let $A$ be a $\Phi$-ring of $X$. Then $(A^{\triangleleft})_x \cong (A_x)^{\triangleleft}$ and $(A^{\mathrm{pur}})_x \cong (A_x)^{\mathrm{pur}}$.
\end{prop}

This follows from the proofs of \ref{chap1sheafpurif} and \ref{chap1sheafclos} since the purification and closure functors on $\Phi$-rings have right adjoints and thus commute with filtered colimits.

\section{Valuative spaces as $\Phi$-localizations \label{chap1valspaces}}

Let $X$ be a locally $\Phi$-ringed topological space (cf. \ref{chap1localphiringed}). 

\subsection{\label{chap10.2.1}} We define the relative valuative spectrum $\widetilde{X}$ of $X$ to be the set of triples $\tilde{x}=(x,\p,R)$, where $x$ is a point of $X$, $\p$ is a prime ideal of $\Ow_{X,x}^{\triangleleft}$ (cf. \ref{chap1sheafclos}, \ref{chap1sheafstalk}) which does not contain an admissible ideal of $\Ow_{X,x}^{\triangleleft}$ and $R$ is a valuation subring of $\ka(\p) = \mathrm{Frac}(\Ow_{X,x}^{\triangleleft}/ \p)$, such that the following two conditions are satisfied:
\begin{itemize}
\item[$(a)$] The homomorphism $\Ow_{X,x} \rightarrow \ka(\p)$ factors through a local homomorphism $\Ow_{X,x} \rightarrow R$.
\item[$(b)$] For any $a \in \ka(\p)$ there is an admissible ideal $I$ of $\Ow_{X,x}$ satisfying $Ia \subseteq R$. 
\end{itemize}

If $\tilde{x}=(x,\p,R)$ is a point of $\widetilde{X}$, we write $\Gamma(\tilde{x}) = \kappa(\p)^{\times}/R^{\times}$, and $\Gamma(\tilde{x})_+ = \Gamma(\tilde{x}) \sqcup \{ 0 \}$. The multiplicative valuation on $\Ow_{X,x}^{\triangleleft}$ associated to $R$ is a multiplicative map $\Ow_{X,x}^{\triangleleft} \rightarrow \Gamma(\tilde{x})_+$ denoted by $f \mapsto |f(\tilde{x})|$. 

\begin{rema} Here and hereafter we use the term ``valuation'' as a synonym for ``higher-rank norm''. One recovers the usual notion of valuation by reversing the order on $\Gamma(\tilde{x})_+$ and by relabelling $0 \in \Gamma(\tilde{x})_+$ as $\infty$.
\end{rema} 

\subsection{\label{chap10.2.2}} Let $\pi : \widetilde{X} \rightarrow X$ be the map $(x,\p,R) \mapsto x$. If $U$ is an open subset of $X$, $I$ is a finitely generated ideal of $\Gamma(U,\Ow_X)$ and $g$ is an element of $\Gamma(U,\Ow_X)$ such that $(I,g)$ is an admissible ideal of $\Gamma(U,\Ow_X)$ then we set
$$
U \left( g^{-1} I \right) = \{ \tilde{x} \in \pi^{-1}(U) \ | \ |I(\tilde{x})| \leq |g(\tilde{x})|\}.
$$
Here, by $|I(\tilde{x})|$ we denote the maximum of $|f(\tilde{x})|$ where $f$ runs over all elements of $I$; this is equal to $\max(|f_i(\tilde{x})|)_{i=1}^r$ whenever $I$ is generated by $f_1,\dots,f_r$. We endow $\widetilde{X}$ with the topology generated by subsets of the form $U \left( g^{-1} I \right)$. 

\begin{prop}\label{chap10.2.2.1} The map $\pi$ is continuous and its image is equal to the support of the sheaf $\Ow_X^{\triangleleft}$, or equivalently of its subsheaf of rings $\Ow_X^{\mathrm{pur}}$. 
\end{prop}

Since $\Ow_X^{\mathrm{pur}}$ is a subsheaf of rings of $\Ow_X^{\triangleleft}$, its support is equal to the support of $\Ow_X^{\triangleleft}$. For any open subset $U$ of $X$, the inverse image $\pi^{-1}(U) = U(1^{-1} \Gamma(U,\Ow_X))$ is open in $\widetilde{X}$, hence the continuity of $\pi$.  Next we notice that if $\Ow_{X,x}^{\triangleleft}$ vanishes then this ring has no prime ideals, so that $\pi^{-1}(x)$ is empty. Let $x$ be a point of $X$ such that $\Ow_{X,x}^{\triangleleft}$ is nonzero, and let $\p'$ be a minimal prime ideal of $\Ow_{X,x}^{\triangleleft}$. Were $\p'$ to contain an admissible ideal $I$ of $\Ow_{X,x}^{\triangleleft}$, this finitely generated ideal would be nilpotent in the localization $(\Ow_{X,x}^{\triangleleft})_{\p'}$ and we would obtain an element $f $ of $\Ow_{X,x}^{\triangleleft} \setminus \p'$ and an integer $N$ such that $f I^N = 0$, which would contradict the fact that $\Ow_{X,x}^{\triangleleft}$ is $1$-deep. Thus $\p'$ does not contain any admissible ideal of $\Ow_{X,x}^{\triangleleft}$.
Let $R'$ be a valuation subring of $\ka(\p')$ dominating the image of $\Ow_{X,x}$. The construction \ref{chap1Constr} applied to the local homomorphism $\Ow_{X,x} \rightarrow R'$ yields a pair $(\p,R)$ such that $(x,\p,R)$ belongs to $\widetilde{X}$, and therefore the fiber $\pi^{-1}(\{x \})$ is not empty.

\subsection{\label{chap10.2.3}} Let us consider the following sheaf of monoids on $\widetilde{X}$:
$$
S : U \mapsto S(U) = \{ s \in \Gamma(U, \pi^{-1}\Ow_X^{\triangleleft}) \ | \ \forall \tilde{x} \in U, \ |s(\tilde{x})| > 0 \}.
$$
We define $\Ow_{\widetilde{X}}^{\triangleleft}$ to be the sheafification of the presheaf of rings $U \mapsto S(U)^{-1} \Gamma(U, \pi^{-1}\Ow_X^{\triangleleft})$. 

\begin{prop}\label{chap10.2.3.1} Let $\tilde{x}=(x,\p,R)$ be a point of $\widetilde{X}$. The stalk $\Ow_{\widetilde{X},\tilde{x}}^{\triangleleft}$ of $\Ow_{\widetilde{X}}^{\triangleleft}$ at $\tilde{x}$ is the localization at $\p$ of $\Ow_{X,x}^{\triangleleft}$.
\end{prop}

The stalk at $\tilde{x}$ of the presheaf given by $U \mapsto S(U)^{-1} \Gamma(U, \pi^{-1}\Ow_X^{\triangleleft})$ is canonically isomorphic to $S_{\tilde{x}}^{-1} \Ow_{X,x}^{\triangleleft}$ where $S_{\tilde{x}}$ is the stalk at $\tilde{x}$ of the sheaf $S$. If $|s(\tilde{x})|>0$ for a section $s$ of $\Ow_{\widetilde{X}}^{\triangleleft}$ on a neighbourhood $U$ of $x$ then up to shrinking $U$ if necessary we have $|s(\tilde{x})| \geq |I(\tilde{x})|$ for some admissible ideal $I$ of $\Gamma(U,\Ow_{X})$ by \ref{chap10.2.1}$(b)$, and thus $s$ is a section of $S$ on the open set $\{ \tilde{y} \in \pi^{-1}(U) \ | \ |I(\tilde{y})| \leq |s(\tilde{y})|\}$. We thus have
$$
S_{\tilde{x}} = \{ s \in \Ow_{X,x}^{\triangleleft} \ | \ |s(\tilde{x})| > 0 \}.
$$
Since $S_{\tilde{x}}$ is the complement of $\p$ we obtain $\Ow_{\widetilde{X},\tilde{x}}^{\triangleleft} \cong (\Ow_{X,x}^{\triangleleft})_{\p}$.

For each point $\widetilde{x}$ of $\widetilde{X}$, the valuation $|\cdot(\tilde{x})|$ on $\Ow_{X,x}^{\triangleleft}$ extends uniquely to its localization $\Ow_{\widetilde{X},\tilde{x}}^{\triangleleft}$. This allows us to consider the subsheaf $\Ow_{\widetilde{X}}$ of $\Ow_{\widetilde{X}}^{\triangleleft}$ defined as follows:
$$
\Ow_{\widetilde{X}} : U \mapsto \Ow_{\widetilde{X}}(U) = \{ f \in \Ow_{\widetilde{X}}^{\triangleleft}(U) \ | \ \forall \tilde{x} \in U, \ |f(\tilde{x})| \leq 1 \}.
$$
We endow $\Ow_{\widetilde{X}}$ with a structure of a presheaf of $\Phi$-rings by declaring a finitely generated $I$ of $\Ow_{\widetilde{X}}(U)$ to be admissible whenever $|I(\tilde{x})| >0$ at each point $\tilde{x}$ of $U$. Then $\Ow_{\widetilde{X}}$ is a sheaf of $\Phi$-rings on $\widetilde{X}$ by \ref{chap1sheafcond}.

\begin{prop}\label{chap1stalkstruct}Let $\tilde{x}=(x,\p,R)$ be a point of $\widetilde{X}$. 
\begin{itemize}
\item[$(i)$] The stalk $\Ow_{\widetilde{X},\tilde{x}}$ of $\Ow_{\widetilde{X}}$ at $\tilde{x}$ is the inverse image of $R \subseteq \ka(\p)$ in $\Ow_{\widetilde{X},\tilde{x}}^{\triangleleft}$.
\item[$(ii)$] The stalk $\Ow_{\widetilde{X},\tilde{x}}$ is a $\Phi$-local $\Phi$-ring with residual valuation ring $R$ (cf. \ref{chap12.1}).
\end{itemize} 
\end{prop}

The stalk $\Ow_{\widetilde{X},\tilde{x}}$ is the preimage of $R$ under the canonical surjective homomorphism $(\Ow_{X,x}^{\triangleleft})_{\p} \rightarrow \ka(\p)$, and a finitely generated ideal $I$ of $\Ow_{\widetilde{X},\tilde{x}}$ is admissible if and only if $|I(\tilde{x})| > 0$, i.e. if and only if $I$ is generated by an element which is invertible in $ (\Ow_{X,x}^{\triangleleft})_{\p}$. By \ref{chap1remarkphilocal} one concludes that $\Ow_{\widetilde{X},\tilde{x}}$ is a $\Phi$-local $\Phi$-ring.

\begin{cor} The canonical morphism $(\Ow_{\widetilde{X}})^{\triangleleft} \rightarrow  \Ow_{\widetilde{X}}^{\triangleleft}$ of $\Phi$-rings is an isomorphism.
\end{cor}

This follows from \ref{chap1sheafstalk} and \ref{chap1stalkstruct}.

\begin{defi}\label{chap1phinormali} The \textbf{$\Phi$-localization} of a locally $\Phi$-ringed topological space $X$ is the $\Phi$-locally $\Phi$-ringed topological space $\widetilde{X}$ endowed with the sheaf of $\Phi$-rings $\Ow_{\widetilde{X}}$.
\end{defi}

\subsection{\label{chap10.2.4}} The canonical morphism $\pi^{-1}\Ow_X \rightarrow \Ow_{\widetilde{X}}^{\triangleleft}$ uniquely factors through a morphism of sheaves of $\Phi$-rings from $\pi^{-1}\Ow_X$ to $\Ow_{\widetilde{X}}$. This provides $\pi$ with the structure of a morphism of locally $\Phi$-ringed topological spaces (cf. \ref{chap1stalkstruct}). Moreover the source $\widetilde{X}$ of $\pi$ is $\Phi$-locally $\Phi$-ringed by \ref{chap1stalkstruct}.

\begin{prop}\label{chap1universalprop} The morphism $\pi : \widetilde{X} \rightarrow X$ is terminal among morphisms of locally $\Phi$-ringed topological spaces from a $\Phi$-locally $\Phi$-ringed topological space to $X$, and hence is an isomorphism if the stalks of $\Ow_X$ are $\Phi$-local. In particular, the assignment $X \mapsto \widetilde{X}$ is functorial and provides a right adjoint to the canonical inclusion functor from $\Phi$-locally $\Phi$-ringed topological spaces to locally $\Phi$-ringed topological spaces.
\end{prop}

Let us consider a morphism $\varphi : Y \rightarrow X$ of locally $\Phi$-ringed topological spaces such that the stalks of $\Ow_Y$ are $\Phi$-local $\Phi$-rings. For each point $y$ of $Y$ we denote by $f \mapsto |f(y)|$ the valuation on $\Ow_{Y,y}^{\triangleleft}$ associated to its residual valuation ring (cf. \ref{chap12.1}).

We first construct a map $\widetilde{\varphi} : Y \rightarrow \widetilde{X}$ such that $\pi \widetilde{\varphi} = \varphi$. Let $y$ be a point of $Y$ and consider the morphism of $\Phi$-rings $\Ow_{X,\varphi(y)} \rightarrow \Ow_{Y,y}$. Let $\mathfrak{m}_y$ be the maximal ideal of $\Ow_{Y,y}^{\triangleleft}$, and let $(\p,R)$ be the pair obtained from the local homomorphism $\Ow_{X,\varphi(y)} \rightarrow \Ow_{Y,y} / \mathfrak{m}_y$ by the construction \ref{chap1Constr}. We set $\widetilde{\varphi}(y) = (\varphi(y), \p, R)$.

The map $\widetilde{\varphi}$ is continuous. Indeed if $I = (f_1,\dots,f_r)$ is an admissible ideal of $\Ow_X(U)$ and if $g$ is an element of $\Ow_X(U)$ then for each $y \in \pi^{-1}(U)$ we have $|\varphi^{-1}I(y)| \leq |\varphi^{-1}g(y)|$ if and only if $|I(\widetilde{\varphi}(y))| \leq |g(\widetilde{\varphi}(y))|$ (cf. \ref{chap1Constr}(d)). We thus have
$$
\widetilde{\varphi}^{-1}\left( U(g^{-1} I) \right) = \{ y \in \varphi^{-1}(U) \ | \ |\varphi^{-1}I(y)| \leq |\varphi^{-1}g(y)| \}. 
$$
This is an open subset of $Y$ since if $y$ is one of its points then $f_i = a_i g$ for some elements $a_1,\dots,a_r$ of $\Ow_{Y,y}$ (cf. \ref{chap12.1}), so that $I \Ow_V \subseteq g \Ow_V$ for some open neighbourhood $V$ of $y$ in $\varphi^{-1}(U)$.

Let us consider the homomorphism
$$
\widetilde{\varphi}^{-1} \pi^{-1} \Ow_X^{\triangleleft}   \xrightarrow{\sim} \varphi^{-1}  \Ow_X^{\triangleleft} \rightarrow \Ow_{Y}^{\triangleleft}.
$$
This homomorphism sends the subsheaf of monoids $\widetilde{\varphi}^{-1} S$ (cf. \ref{chap10.2.3}) into the sheaf of invertible elements of $\Ow_{Y,y}^{\triangleleft}$ (cf. \ref{chap12.1}), so that it uniquely factors though a local homomorphism $\widetilde{\varphi}^{-1} \Ow_{\widetilde{X}}^{\triangleleft} \rightarrow \Ow_{Y}^{\triangleleft}$. Moreovever the image of $\widetilde{\varphi}^{-1} \Ow_{\widetilde{X}}$ is contained in the subsheaf $\Ow_{Y,y}$, so that we obtain a morphism $\widetilde{\varphi}^{-1} \Ow_{\widetilde{X}} \rightarrow \Ow_{Y}$ which endows $\widetilde{\varphi}$ with a structure of morphism of locally $\Phi$-ringed topological spaces.
%
%\begin{cor} If $X$ is $\Phi$-locally $\Phi$-ringed then $\pi$ is an isomorphism of locally $\Phi$-ringed topological spaces.
%\end{cor}
%
%This follows from the universal property \ref{chap1universalprop}.
%
%\begin{cor}\label{chap1fonct} If $Y \rightarrow X$ is a morphism of locally $\Phi$-ringed topological spaces, then there is a unique morphism $\widetilde{Y} \rightarrow \widetilde{X}$ of locally $\Phi$-ringed topological spaces such that the diagram 
%\begin{equation*}\label{chap1intro3b}
%\xymatrix{\widetilde{Y}\ar[r]\ar[d]&Y\ar[d]\\
%\widetilde{X}\ar[r]&X}
%\end{equation*}
%is commutative.
%\end{cor} 
%
%This follows from from the universal property \ref{chap1universalprop} by considering the composed morphism $\widetilde{Y} \rightarrow Y \rightarrow X$. 
%
%As a direct consequence of \ref{chap1universalprop} and \ref{chap1fonct}, we have the following result:
%
%\begin{cor}\label{chap1fonct2} The canonical inclusion functor from $\Phi$-locally $\Phi$-ringed topological spaces to locally $\Phi$-ringed topological spaces has a right adjoint, the \textbf{$\Phi$-localization functor}, given on objects by $X \mapsto \widetilde{X}$.
%\end{cor} 

\subsection{\label{chap10.2.5}} We assume for the remainder of this section that the locally ringed topological space underlying $X$ is a scheme. We make the following two definitions:

 \begin{itemize}
\item[$\triangleright$] A quasi-coherent ideal $\mathcal{I} \subseteq \Ow_{X}$ is \textbf{admissible} if for any affine open subscheme $V \subseteq X$, the ideal $\mathcal{I}(V)$ is admissible in the $\Phi$-ring $\Ow_X(V)$.
\item[$\triangleright$] An \textbf{admissible blow-up} $X' \rightarrow X$ is the blow-up of $X$ along a closed subscheme defined by an admissible quasi-coherent ideal sheaf. 
 \end{itemize}

Let $\pi : \widetilde{X} \rightarrow X$ be the $\Phi$-localization of $X$ (cf. \ref{chap1phinormali}). Let $\mathcal{I}$ be an admissible quasi-coherent ideal sheaf on $X$, and let $V$ be an affine open subset of $X$. Then $\mathcal{I}(V) = (f_1,\dots,f_r)$ is an admissible ideal in $\Ow_X(V)$. Hence $\mathcal{I}(V)$ becomes invertible on $\pi^{-1}(V)$. More precisely, we have the decomposition by \ref{chap12.1}
$$
\pi^{-1}(V) = \bigcup_{i=1}^r V(f_i^{-1} \mathcal{I}(V)),
$$
and the restriction of $\mathcal{I} \Ow_{\widetilde{X}}$ to $V(f_i^{-1} \mathcal{I}(V))$ is freely generated by $f_i$. Thus $\mathcal{I} \Ow_{\widetilde{X}}$ is a locally free  $\Ow_{\widetilde{X}}$-module of rank $1$. In particular if we denote by $X_{\mathcal{I}}$ the blow-up in $X$ of the closed subscheme defined by $\mathcal{I}$, then from the universal property of the blow-up we obtain a unique factorisation $\widetilde{X} \rightarrow X_{\mathcal{I}} \rightarrow X$ of $\pi$.\footnote{The universal property of the blow-up, while usually stated in the category of schemes, also holds in the larger category of locally ringed topological spaces} This yields a morphism
\begin{equation} \label{chap1comparisonmap}
\widetilde{X} \rightarrow \lim_{X' \rightarrow X} \ X'
\end{equation}
of locally $\Phi$-ringed topological spaces, where the limit runs over the category of admissible blow-ups of $X$. The latter category is cofiltered; indeed it is non empty since the identity from $X$ to itself is an admissible blow-up, any pair $X_{\mathcal{I}},X_{\mathcal{J}}$ of admissible blow-ups of $X$ is dominated by the admissible blow-up $X_{\mathcal{I} \mathcal{J}}$ and there are no pair of distinct parallel arrows since for any morphism of schemes $f : X' \rightarrow X$, the universal property of the blow-up ensures that there exists a morphism of $X$-schemes from $X'$ to $X_{\mathcal{I}}$ if and only if $f^* \mathcal{I}$ is invertible, in which case there exists a unique such morphism.

The limit (\ref{chap1comparisonmap}) of locally $\Phi$-ringed topological spaces is given by taking the corresponding limit in the category of locally ringed topological spaces, and by declaring the canonical isomorphism of its structure sheaf with the colimit
$$
\mathrm{colim}_{f : X' \rightarrow X} \ \ f^{-1} \Ow_{X'}
$$
to be an isomorphism of sheaves of $\Phi$-rings.

\begin{teo}\label{chap1comparisoneclate} Assume that $X$ is quasi-compact and quasi-separated and that there exists a quasi-compact open subscheme $U \subseteq X$ such that for any open subset $V \subseteq X$, a finitely generated ideal $I$ of $\Ow_X(V)$ is admissible if and only if $I \Ow_{U \cap V} = \Ow_{U \cap V}$. Then the canonical morphism
$$
\widetilde{X} \rightarrow \lim_{X' \rightarrow X} \ X'
$$
from $\widetilde{X}$ to the inverse limit of the admissible blow-ups of $X$ is an isomorphism of locally $\Phi$-ringed topological spaces. In particular $\widetilde{X}$ is quasi-compact and is even a spectral space, cf. $($\cite{Stacks} $\mathrm{0A2V}$, $\mathrm{0A2Z})$
\end{teo}

For any admissible blow-up $f : X' \rightarrow X$ we endow $X'$ with the unique structure of locally $\Phi$-ringed topological space such that for any open subset $V \subseteq X'$, a finitely generated ideal $I$ of $\Ow_{X'}(V)$ is admissible if and only if $I \Ow_{f^{-1}(U) \cap V} = \Ow_{f^{-1}(U) \cap V}$.

In order to prove \ref{chap1comparisoneclate}, it is sufficient to show that the limit in (\ref{chap1comparisonmap}), which is a locally $\Phi$-ringed topological space, is actually $\Phi$-locally $\Phi$-ringed. Indeed, a morphism from a $\Phi$-locally $\Phi$-ringed topological space to $X$ uniquely factors through the limit in (\ref{chap1comparisonmap}) by the universal property of the blow-up, hence if the latter is $\Phi$-locally $\Phi$-ringed then it is terminal among morphisms of locally $\Phi$-ringed topological spaces from a $\Phi$-locally $\Phi$-ringed topological space to $X$, so that the canonical morphism (\ref{chap1comparisonmap}) is an isomorphism by virtue of \ref{chap1universalprop}.

Let us therefore prove that the limit in (\ref{chap1comparisonmap}) is $\Phi$-locally $\Phi$-ringed. Let $x = (x_{\mathcal{I}})_{\mathcal{I}}$ be a point on this limit. The stalk of the structure sheaf at $x$ is the colimit of $\mathcal{I} \mapsto  \Ow_{ X_{\mathcal{I}} , x_{\mathcal{I}} }$ where the index $\mathcal{I}$ runs over all admissible quasi-coherent ideal sheaves on $X$. Let $J$ be an admissible ideal of $\Ow_{ X_{\mathcal{I}} , x_{\mathcal{I}} }$. We must show that $J \Ow_{ X_{\mathcal{IK}} , x_{\mathcal{IK}} }$ is an invertible ideal for some admissible quasi-coherent ideal sheaf $\mathcal{K}$.

We first prove that there exists an admissible quasi-coherent ideal sheaf $\mathcal{J}$ on $X_{\mathcal{I}}$ such that $\mathcal{J}_{x_{\mathcal{I}}} = J$. The ideal $J$ is the stalk at $x_{\mathcal{I}}$ of an admissible quasi-coherent ideal sheaf $\mathcal{J}_0$ on an open neighbourhood $V$ of $x_{\mathcal{I}}$ in $X_{\mathcal{I}}$. Since $\mathcal{J}_0 \Ow_{U \cap V} = \Ow_{U \cap V}$ where $U$ denotes the (isomorphic) inverse image of $U$ in $X_{\mathcal{I}}$, there exists a unique quasi-coherent ideal sheaf $\mathcal{J}_1$ on $U \cup V$ such that $\mathcal{J}_{1|U} = \Ow_U$ and $\mathcal{J}_{1|V} = \mathcal{J}_0$. Since $X_{\mathcal{I}}$ and $U$ are quasi-compact, there exists by (\cite{EGA1}, $9.4.7$) a quasi-coherent ideal sheaf $\mathcal{J}$ of finite type on $X_{\mathcal{I}}$ such that $\mathcal{J}_{|U \cup V} = \mathcal{J}_{1}$. The condition $\mathcal{J} \Ow_{U} = \Ow_U$ implies that $\mathcal{J}$ is admissible.

By (\cite{Raynaud-Gruson71}, $I.5.1.4$) the composition of blow-ups $(X_{\mathcal{I}} )_{\mathcal{J}} \rightarrow X$ is isomorphic to the blow-up of an admissible quasi-coherent ideal sheaf $\mathcal{K}$ on $X$. Thus $J \Ow_{ X_{\mathcal{IK}} , x_{\mathcal{IK}} }$ is an invertible ideal and we have shown that the stalks of the structure sheaf of the limit appearing in (\ref{chap1comparisonmap}) are $\Phi$-local $\Phi$-rings. As explained above, this together with \ref{chap1universalprop} and the universal property of the blow-up concludes the proof of Theorem \ref{chap1comparisoneclate}.

\section{The flattening property of $\Phi$-localizations \label{chap1flattprop}}

\subsection{\label{chap10.3.1}} If $M$ is a module over a valuation ring $R$ then $M$ is $R$-flat if and only if for any nonzero element $r$ of $R$ the module $M$ has no nonzero $r$-torsion, see (\cite{Bourbaki}, VI \S 3.6 Lemme 1) or (\cite{Stacks} $\mathrm{0539}$). We generalize this fact to modules over $\Phi$-local $\Phi$-rings:

\begin{prop}\label{chap10.3.2} Let $A$ be a $\Phi$-local $\Phi$-ring, and let $M$ be a $A$-module. The following are equivalent:
\begin{itemize}
\item[$(i)$] The $A$-module $M$ is flat.
\item[$(ii)$] The $A^{\triangleleft}$-module $M^{\triangleleft} = M \otimes_{A} A^{\triangleleft}$ is flat, and the map $M \rightarrow M^{\triangleleft}$ is injective.
\end{itemize}
\end{prop}

We have $(i) \implies (ii)$, since $A$ is $1$-deep by \ref{chap12.1}. We thus focus on the converse implication. We use the equational criterion of flatness, see (\cite{Bourbaki}, I \S 2.11 Corollaire 1) or (\cite{Stacks} $\mathrm{00HK}$). Namely, given a relation $\sum_{d \in D} a_d x_d = 0$ with $a_d$ in $A$ and $x_d$ in $M$, for some finite set $D$, we would like to show that this relation is trivial, in the sense that we can find relations of the form
$$
x_d = \sum_{e \in E} b_{de} y_e,
$$
such that $\sum_{d \in D} a_d b_{de} = 0$ for any $e$. We proceed by induction on the cardinality of $D$, the case $D = \emptyset$ being empty. 

First, assume that we have a relation of the form $\sum_{d \in D} a_d c_d =0$, where $(c_d)_{d \in D} \subseteq A$ and $c_{d_0} =1$ for some $d_0$ in $D$. We then have
$$
\sum_{d \in D \setminus \{ d_0\}} a_d (x_d - c_d x_{d_0}) = \sum_{d \in D } a_d (x_d - c_d x_{d_0}) = 0.
$$
Since this relation has fewer terms than the original one, the induction hypothesis ensures that we can find relations
$$
x_d - c_d x_{d_0}= \sum_{e \in E} b_{de} y_e
$$
with $\sum_{d \in D \setminus \{ d_0\}} a_d b_{de} = 0$ for any $e$. We set $E'=E \sqcup \{ \star \}$, $y_{\star} = x_{d_0}$, $b_{d_0 e} = 0$ for all $e \in E$, and $b_{d \star} = c_d$ for all $d \in D$, so that $b_{d_0 \star} = 1$. We then have 
$$
x_d = \sum_{e \in E'} b_{de} y_e,
$$
for any $d$ in $D$, and $\sum_{d \in D} a_d b_{de} = 0$ for any $e$ in $E'$. Thus our relation is trivial.

Let us now assume that the previous case does not happen. Since $M^{\triangleleft}$ is flat over $A^{\triangleleft}$, our relation is trivial in $M^{\triangleleft}$, hence we can find relations
$$
x_d = \sum_{e \in E} b_{de} y_e,
$$
such that $\sum_{d \in D} a_d b_{de} = 0$ for any $e$, with $y_e$ in $M$. Here the coefficients $b_{de}$ are elements of $A^{\triangleleft}$, hence it is sufficient for our purpose to show that they belong to $A$.

Actually, we have the stronger result that for all $(d,e) \in D \times E$, the element $b_{de}$ belongs to the maximal ideal of $A^{\triangleleft}$. Indeed, if it were not the case, then, for some $e_0 \in E$, the ideal generated by $(b_{de_0})_{d \in D}$ would not be contained in the maximal ideal of $A^{\triangleleft}$. We would then have an element $s$ of $A$, invertible in $A^{\triangleleft}$, such that all the $sb_{de_0}$ belong to $A$. The ideal of $A$ generated by the family $(sb_{de_0})_{d \in D}$ would then be admissible, hence generated by $sb_{d_0e_0}$ for some $d_0$ in $D$. This would imply that $sb_{de_0} = c_d sb_{d_0 e_0}$ for some $c_d$ in $A$, with $c_{d_0} = 1$, and the relation $\sum_{d \in D} a_d b_{de_0} = 0$ would yield $\sum_{d \in D} a_d c_d = 0$. This is a contradiction.

The following lemma is both a consequence and a generalization of (\cite{Stacks} $\mathrm{053E}$):

\begin{prop}\label{chap12.2.5} Let $A$ be a $\Phi$-local $\Phi$-ring, let $A \rightarrow B$ be a ring homomorphism of finite type, and let $M$ be a $B$-module of finite type. If $M$ is flat over $A$, and if $M \otimes_{A} A^{\triangleleft}$ is a finitely presented $B \otimes_{A} A^{\triangleleft}$-module, then $M$ is a finitely presented $B$-module.
\end{prop}

Since $M$ is a $B$-module of finite type, we can find an exact sequence
$$
0 \rightarrow K \rightarrow B^{\oplus r} \rightarrow M \rightarrow 0.
$$
Let $\m$ be the maximal ideal of $A^{\triangleleft}$, which is contained in $A$ by \ref{chap12.1} (ii), and let $S = A \setminus \m$. Since $M$ is flat over $A$, the following exact sequence is still exact:
$$
0 \rightarrow K/\m K \rightarrow \left(B/\m B\right)^{\oplus r} \rightarrow M/ \m M \rightarrow 0.
$$
By (\cite{Stacks} $\mathrm{053E}$), the kernel $K/ \m K$ is a finitely generated $B/ \m B$-module. Thus there is a finite subset $\mathfrak{S}$ of $K$ such that $K = B \mathfrak{S} + \m K $. But we also have an exact sequence
$$
0 \rightarrow S^{-1} K \rightarrow \left(S^{-1}B \right)^{\oplus r} \rightarrow S^{-1} M  \rightarrow 0,
$$
and since by hypothesis $S^{-1} M = M \otimes_{A} A^{\triangleleft}$ is a finitely presented module over $S^{-1} B = B \otimes_{A} A^{\triangleleft}$, the kernel $S^{-1} K $ is a finitely generated $S^{-1} B$-module. Thus, by enlarging $\mathfrak{S}$ if necessary, we can also assume that $S^{-1} K = S^{-1}B \mathfrak{S} $. We claim that $K = B \mathfrak{S}$. Indeed, the relations $S^{-1} K = S^{-1}B \mathfrak{S} $ and $\m = S^{-1} \m$ imply
$$
\m K = \m S^{-1} K = \m  S^{-1}B \mathfrak{S} = \m  B \mathfrak{S},
$$
so that
$$
K =B \mathfrak{S} + \m K =  B \mathfrak{S} + \m B \mathfrak{S} =  B \mathfrak{S}.
$$
In particular, $K$ is a finitely generated $B$-module, and consequently $M$ is a finitely presented $B$-module.

\begin{cor}\label{chap12.3} Let $A$ be a $\Phi$-local $\Phi$-ring, let $A \rightarrow B$ be a ring homomorphism of finite type. If $B$ is flat over $A$, and if $B \otimes_{A} A^{\triangleleft}$ is a finitely presented $A^{\triangleleft}$-algebra, then $B$ is a finitely presented $A$-algebra.
\end{cor}

Indeed, if we choose a surjection $B' = A[X_1,\dots,X_n]\rightarrow B$, then $B$ is flat over $A$, and $B \otimes_{A} A^{\triangleleft}$ is a finitely presented $B' \otimes_{A} A^{\triangleleft}$-module since $B \otimes_{A} A^{\triangleleft}$ is a finitely presented $A^{\triangleleft}$-algebra. Thus, by \ref{chap12.2.5}, we conclude that $B$ is a finitely presented $B'$-module, and thus that $B$ is a finitely presented $A$-algebra.

\subsection{\label{chap10.3.7}} Let $(X,\Ow_X)$ be a ringed topological space. A commutative $\Ow_X$-algebra $\mathcal{A}$ is said to be \textit{finitely presented} if for any point $x$ of $X$ there exists an open neighbourhood $U$ of $x$ in $X$ and an isomorphism
$$
A_U \otimes_{\Gamma(U,\Ow_X)} \Ow_U \rightarrow \mathcal{A}_{|U},
$$
of $\Ow_U$-algebras, for some finitely presented $\Gamma(U,\Ow_X)$-algebra $A_U$.

\begin{prop}\label{chap1stacksargument} Let $(X,\Ow_X)$ be a ringed topological space, let $\mathcal{A}$ be finitely presented commutative $\Ow_X$-algebra, and let $\mathcal{G}$ be a finitely presented $\mathcal{A}$-module. Let $x$ be a point of $X$ such that the stalk $\mathcal{M}_x$ is a flat $\Ow_{X,x}$-module. Then there exists an open neighbourhood $U$ of $x$ such that the restriction $\mathcal{M}_{|U}$ is a flat $\Ow_U$-module.
\end{prop}

Up to replacing $X$ by some open neighbourhood of $x$ in $X$, we can assume (and we do) that $\mathcal{A}$ is of the form $A \otimes_{\Gamma(X,\Ow_X)} \Ow_X$, for some finitely presented $\Gamma(X,\Ow_X)$-algebra $A$, and that we have a global presentation
$$
\mathcal{A}^{\oplus m} \xrightarrow[]{\varphi} \mathcal{A}^{\oplus n}   \rightarrow \mathcal{M}  \rightarrow 0 .
$$
For any integer $j$ between $1$ and $m$, let $(a_{i,j})_{i=1}^n$ be the image by $\varphi$ of the $j$-th basis vector. Each $a_{i,j}$ belongs to $\Gamma(X,\mathcal{A})$. Let $V$ be an open neighbourhood of $x$ in $X$ such that each restriction $a_{i,j|V}$ is the image in $\Gamma(V,\mathcal{A})$ of an element $b_{i,j}$ of the finitely presented $\Gamma(V,\Ow_X)$-algebra $A_V = A  \otimes_{\Gamma(X,\Ow_X)} \Gamma(V, \Ow_X)$. Let $M_V$ be the cokernel of the matrix $(b_{i,j})_{i,j}$. It is a finitely presented $A_V$-module, and we have an isomorphism
$$
M_V \otimes_{A_V} \mathcal{A}_{|V} \rightarrow \mathcal{M}_{|V},
$$
of $\mathcal{A}_{|V}$-modules. For each open neighbourhood $U$ of $x$ in $V$, the $\Gamma(U,\Ow_X)$-algebra $A_U = A_V  \otimes_{\Gamma(V,\Ow_X)} \Gamma(U, \Ow_X)$ is finitely presented, and the $A_U$-module $M_U = M_V \otimes_{A_V} A_U$ is finitely presented. Moreover, by hypothesis the filtered colimit
$$
\mathrm{colim}_{U} M_U \cong \mathcal{M}_x,
$$
where $U$ runs over the cofiltered set of open neighbourhood of $x$ in $V$, is a flat module over the colimit 
$$
\mathrm{colim}_U \Gamma(U, \Ow_X)  = \Ow_{X,x}.
$$
By (\cite{Stacks} 02JO(3)), this implies that there exists an open neighbourhood $U$ of $x$ in $V$ such that $M_U$ is a flat $\Gamma(U, \Ow_X)$-module. In particular, the $\Ow_U$-module 
$$
\mathcal{M}_{|U} \cong M_U \otimes_{A_U}  \mathcal{A}_{|U} \cong M_U \otimes_{\Gamma(U,\Ow_X)} \Ow_U,
$$
is flat, hence the conclusion of Proposition \ref{chap1stacksargument}.

\subsection{\label{chap10.3.3}} Let $X$ be a $\Phi$-ringed topological space, and let $\mathcal{F}$ be an $\Ow_X$-module. We define the \textbf{closure} $\mathcal{F}^{\triangleleft}$ of $\mathcal{F}$ to be the tensor product $\mathcal{F} \otimes_{\Ow_X} \Ow_{X}^{\triangleleft}$ and the \textbf{purification} $\mathcal{F}^{\pur}$ of $\mathcal{F}$ to be the subsheaf of $\mathcal{F}^{\triangleleft}$ generated by the image of $\mathcal{F}$.

\begin{teo}\label{chap1RG1} Let $X$ be a locally $\Phi$-ringed topological space with the $\Phi$-localization $\pi : \widetilde{X} \rightarrow X$ (cf. \ref{chap1phinormali}), and let $\mathcal{F}$ be an $\Ow_X$-module such that $\mathcal{F}^{\triangleleft}$ is a flat $\Ow_X^{\triangleleft}$-module (cf. \ref{chap10.3.3}). Then $(\pi^{*} \mathcal{F})^{\pur}$ is a flat $\Ow_{\widetilde{X}}$-module.
\end{teo}

The stalk of $(\pi^{*} \mathcal{F})^{\pur}$ at a point $\tilde{x} = (x,\p,R)$ of $\widetilde{X}$ is the purification of $\mathcal{F}_x \otimes_{\Ow_{X,x}} \Ow_{\widetilde{X},\tilde{x}}$. Since $(\pi^{*} \mathcal{F})_{\tilde{x}}^{\triangleleft} \cong (\mathcal{F}_x^{\triangleleft})_{\p}$ is a flat module over $\Ow_{\widetilde{X},\tilde{x}}^{\triangleleft} \cong (\Ow_{X,x}^{\triangleleft})_{\p}$ (cf. \ref{chap1stalkstruct}), one concludes by \ref{chap10.3.2} that $(\pi^{*} \mathcal{F})_{\tilde{x}}^{\pur}$ is flat over the $\Phi$-local $\Phi$-ring $\Ow_{\widetilde{X},\tilde{x}}$.

\begin{teo}\label{chap1RG2} Let $X$ be a locally $\Phi$-ringed topological space with the $\Phi$-localization $\pi : \widetilde{X} \rightarrow X$ (cf. \ref{chap1phinormali}). Let $\mathcal{A}$ be a finitely presented commutative $\Ow_X$-algebra (cf. \ref{chap10.3.7}) and let $\mathcal{F}$ be an $\mathcal{A}$-module of finite type such that $\mathcal{F}^{\triangleleft}$ is finitely presented as an $\mathcal{A}^{\triangleleft}$-module and flat as an $\Ow_X^{\triangleleft}$-module (cf. \ref{chap10.3.3}). Then $(\pi^{*} \mathcal{F})^{\pur}$ is finitely presented as a $\pi^{*}\mathcal{A}$-module and flat as an $\Ow_{\widetilde{X}}$-module.
\end{teo}

By using \ref{chap1RG1} and replacing $X$ with $\widetilde{X}$ it is sufficient for the purpose of proving \ref{chap1RG2} to show the following:

\begin{prop}\label{chap1RG3} Let $X$ be a $\Phi$-locally $\Phi$-ringed topological space. Let $\mathcal{A}$ be a finitely presented commutative $\Ow_X$-algebra (cf. \ref{chap10.3.7}) and let $\mathcal{F}$ be an $\mathcal{A}$-module of finite type such that $\mathcal{F}^{\triangleleft}$ is finitely presented as an $\mathcal{A}^{\triangleleft}$-module and $\mathcal{F}$ is flat as an $\Ow_X$-module. Then $\mathcal{F}$ is finitely presented as an $\mathcal{A}$-module.
\end{prop}

Let $x$ be a point of $X$. Let $U$ be an open neighbourhood of $x$ such that there exists a surjective homomorphism $\psi : \mathcal{A}_{|U}^{\oplus n} \rightarrow \mathcal{F}_{|U}$. The $\mathcal{A}^{\triangleleft}_{|U}$-module $\mathcal{F}^{\triangleleft}_{|U}$ is finitely presented, and by (\cite{Stacks} 01BP) this implies that the $\mathcal{A}^{\triangleleft}_{|U}$-module $(\ker \psi)^{\triangleleft}$ is of finite type. By \ref{chap12.2.5}, the $\mathcal{A}_x$-module $\mathcal{F}_x$ is finitely presented, and hence the stalk of $\ker \psi$ at $x$ is a finitely generated $\mathcal{A}_x$-module. Thus up to shrinking $U$ if necessary there exists a finitely generated sub-module $\mathcal{H} \subseteq \ker \psi$ such that $\mathcal{H}_x = (\ker \psi)_x$ and $\mathcal{H}^{\triangleleft} = (\ker \psi)^{\triangleleft}$. Let $\mathcal{G}$ be the quotient of $\mathcal{A}_{|U}^{\oplus n} $ by $\mathcal{H}$. Then $\mathcal{G}$ is a finitely presented $\mathcal{A}_{|U}$-module and we have a surjective homomorphism $\mathcal{G} \rightarrow \mathcal{F}$ such that $\mathcal{G}_x \cong \mathcal{F}_x$ and $\mathcal{G}^{\triangleleft} \cong \mathcal{F}^{\triangleleft}$. Moreover, $\mathcal{G}$ and $\mathcal{F}$ have the same image in $\mathcal{F}^{\triangleleft}$, so that $\mathcal{F}$ is the purification of $\mathcal{G}$. Since $\mathcal{A}_{|U}$ is a finitely presented $\Ow_U$-algebra, since $\mathcal{G}_x$ is a flat $\Ow_{X,x}$-module and since $\mathcal{G}$ is a finitely presented $\mathcal{A}_{|U}$-module, we conclude by Proposition \ref{chap1stacksargument} that there exists an open neighbourhood $V \subseteq U$ of $x$ such that $\mathcal{G}_{|V}$ is a flat $\Ow_V$-module. In particular $\mathcal{G}_{|V}$ is isomorphic to its purification $\mathcal{F}_{|V}$, so that $\mathcal{F}_{|V}$ is finitely presented as an $\mathcal{A}_{|V}$-module.
 
\section{Proof of Raynaud-Gruson's theorem \label{chap1proofofrg}}

\subsection{\label{chap1part1}} We first prove the following variant of Raynaud-Gruson's theorem \ref{chap1RGtheo}:

\begin{teo}\label{chap1RGtheo2} Let $X$ be quasi-compact and quasi-separated scheme and let $U$ be a quasi-compact open subset of $X$. Let $\mathcal{A}$ be a (quasi-coherent) finitely presented commutative $\Ow_X$-algebra and let $\mathcal{F}$ be a quasi-coherent $\mathcal{A}$-module of finite type. Assume that $\mathcal{F}_{|U}$ is finitely presented over $\mathcal{A}_{|U}$ and flat over $\Ow_U$. Then there exists a blow-up $f : X' \rightarrow X$ such that:
\begin{itemize}
\item[$(1)$] The center of $f$ is a finitely presented closed subscheme of $X$, disjoint from $U$.
\item[$(2)$] The strict transform $\mathcal{F}'$ of $\mathcal{F}$ along $f$ is finitely presented over $f^* \mathcal{A}$ and flat over $\Ow_{X'}$.
\end{itemize}
\end{teo}
 
Since $X$ is quasi-compact and quasi-separated and $U$ is quasi-compact, the complement of $U$ is the closed subscheme defined by a finitely generated quasi-coherent ideal sheaf $\mathcal{J}$ on $X$: indeed, the quasi-coherent ideal defining this complement with its reduced structure is a filtered union of the family $(\mathcal{J}_{\lambda})_{\lambda \in \Lambda}$ of its finitely generated quasi-coherent subsheaves by (\cite{EGA1}, $6.9.9$), hence the quasi-compact open subset $U$ is the filtered union of the open subsets $(X \setminus V(\mathcal{J}_{\lambda}))_{\lambda \in \Lambda}$, and thus $U$ is equal to $X \setminus V(\mathcal{J}_{\lambda})$ for some $\lambda$. Let $f : X' \rightarrow X$ be the blowing-up of $X$ along $\mathcal{J}$. For any blow-up $g : X'' \rightarrow X'$ whose center is a finitely presented closed subscheme of $X'$ disjoint from $f^{-1}(U)$, the composition $fg$ is a blow-up whose center is a finitely presented closed subscheme of $X$ disjoint from $U$, cf. (\cite{Raynaud-Gruson71}, $I.5.1.4$). Hence by replacing $X$ with $X'$, we can assume (and we do) that $\mathcal{J}$ is invertible.

We endow $X$ with the structure of a locally $\Phi$-ringed topological space by declaring for any open subset $V \subseteq X$ a finitely generated ideal $I$ of $\Ow_X(V)$ to be admissible whenever $I \Ow_{U \cap V} = \Ow_{U \cap V}$. Since $\mathcal{J}$ is invertible we have $\Ow_{X}^{\triangleleft} = \mathrm{colim}_{N>0}  \ \mathcal{J}^{-N} = j_* \Ow_U$ where $j : U \rightarrow X$ is the canonical inclusion. In particular the assumptions of Theorem \ref{chap1RGtheo2} imply that $\mathcal{F}^{\triangleleft} = \mathrm{colim}_{N>0}  \ \mathcal{F} \otimes_{\Ow_X} \mathcal{J}^{-N}$ is flat over $\Ow_X^{\triangleleft} = \mathrm{colim}_{N>0}  \mathcal{J}^{-N}$ and finitely presented over $\mathcal{A}^{\triangleleft} =  \mathrm{colim}_{N>0}  \ \mathcal{A} \otimes_{\Ow_X} \mathcal{J}^{-N}$.

Let $\pi : \widetilde{X} \rightarrow X$ be the $\Phi$-localization of $X$. By (\cite{EGA1}, $6.9.10$) applied to the relative spectrum of $\mathcal{A}$, there exists a finitely presented $\mathcal{A}$-module $\mathcal{G}$ together with a surjective homomorphism $\mathcal{G} \rightarrow \mathcal{F}$. By Theorem \ref{chap1RG2} the $\pi^{*}\mathcal{A}$-module $(\pi^{*} \mathcal{F})^{\pur}$ is finitely presented, hence the kernel $\mathcal{K}$ of the surjective homomorphism $\pi^* \mathcal{G} \rightarrow (\pi^{*} \mathcal{F})^{\pur}$ is of finite type. 

%
%
%there exists a filtered system $(\mathcal{G}_i)_{i \in I}$ of finitely presented $\mathcal{A}$-modules together with an isomorphism
%$$
%\mathrm{colim} \ \mathcal{G}_i \rightarrow \mathcal{F}.
%$$
%We infer from \ref{chap1RG2} that $(\pi^{*} \mathcal{F})^{\pur}$ is finitely presented as a $\pi^{*}\mathcal{A}$-module (and flat) as an $\Ow_{\widetilde{X}}$-module. Moreover, the $\mathcal{A}^{\triangleleft}$-module $\mathcal{F}^{\triangleleft}$ is finitely presented. By (\cite{Stacks} 01BS), the canonical homomorphisms
%\begin{align*}
%\mathrm{colim} \ \mathrm{Hom}_{\pi^{*}\mathcal{A}}( (\pi^{*} \mathcal{F})^{\pur}, (\pi^{*} \mathcal{G}_i)^{\pur} ) &\rightarrow \mathrm{Hom}_{\pi^{*}\mathcal{A}}( (\pi^{*} \mathcal{F})^{\pur}, (\pi^{*} \mathcal{F})^{\pur}), \\
%\mathrm{colim}  \ \mathrm{Hom}_{\mathcal{A}^{\triangleleft}}( \mathcal{F}^{\triangleleft}, \mathcal{G}_i^{\triangleleft}) &\rightarrow \mathrm{Hom}_{\mathcal{A}^{\triangleleft}}( \mathcal{F}^{\triangleleft}, \mathcal{F}^{\triangleleft}),
%\end{align*}
%are isomorphisms. In particular, there exists a finitely presented $\Ow_X$-module $\mathcal{G}$ together with a morphism $\mathcal{G} \rightarrow \mathcal{F}$ such that $\mathcal{G}^{\triangleleft} \rightarrow \mathcal{F}^{\triangleleft}$ and $(\pi^{*} \mathcal{G})^{\pur} \rightarrow (\pi^{*} \mathcal{F})^{\pur}$ are both split surjections. Since $(\pi^{*} \mathcal{F})^{\pur}$ is a finitely presented $\pi^{*}\mathcal{A}$-module, the kernel $\mathcal{K}$ of the surjective homomorphism $\pi^* \mathcal{G} \rightarrow (\pi^{*} \mathcal{F})^{\pur}$ is of finite type. 

By \ref{chap1comparisoneclate} we have an isomorphism of locally $\Phi$-ringed topological spaces
$$
\widetilde{X} \rightarrow \lim_{X' \rightarrow X} \ X'
$$
from $\widetilde{X}$ to the inverse limit of the admissible blow-ups of $X$ (cf. \ref{chap10.2.5}). We endow each blow-up $g : X' \rightarrow X$ with the structure of a locally $\Phi$-ringed topological space by declaring for any open subset $V \subseteq X'$ a finitely generated ideal $I$ of $\Ow_{X'}(V)$ to be admissible whenever $I \Ow_{g^{-1}(U) \cap V} = \Ow_{g^{-1}(U) \cap V}$. Since $\mathcal{K}$ is a $\pi^{*} \mathcal{A}$-module of finite type, there exists a factorization $\widetilde{X} \xrightarrow[]{\pi'} X' \xrightarrow[]{f} X$ of $\pi$ through an admissible blow-up and a finitely generated sub-$f^* \mathcal{A}$-module $\mathcal{K}'$ of $\mathrm{Ker}(f^* \mathcal{G} \rightarrow (f^{*} \mathcal{F})^{\pur})$ such that the image of $\pi'^* \mathcal{K}'$ in $\pi^* \mathcal{G}$ is equal to $\mathcal{K}$. 

Let us consider the finitely presented $f^* \mathcal{A}$-module $\mathcal{G}' = f^* \mathcal{G} / \mathcal{K}'$. We have $\pi'^* \mathcal{G}' \cong (\pi^{*} \mathcal{F})^{\pur}$, hence an isomorphism $(\pi'^* \mathcal{G}')^{\triangleleft} \cong (\pi^{*} \mathcal{F})^{\triangleleft}$. By \ref{chap10.2.2.1} and \ref{chap10.2.3.1}, the morphism of ringed space 
$$
(\widetilde{X},\Ow_{\widetilde{X}}^{\triangleleft}) \rightarrow (X', \Ow_{X'}^{\triangleleft})
$$
is flat and surjective, hence we deduce that the canonical homomorphism $\mathcal{G}'^{\triangleleft} \rightarrow (f^* \mathcal{F})^{\triangleleft}$ is an isomorphism. Since the homomorphism $\mathcal{G'} \rightarrow f^* \mathcal{F}$ is surjective, we obtain that the purification of $\mathcal{G'}$ is isomorphic to $(f^{*} \mathcal{F})^{\pur}$. More generally, for any admissible blow-up $X''\xrightarrow[]{f'} X'$, the purification of $f'^* \mathcal{G}'$ is isomorphic to $(f'^{*} f^* \mathcal{F})^{\pur}$.

 For each admissible blow-up $X''\xrightarrow[]{f'} X' \xrightarrow[]{f} X$, the set of points $x''$ of $X''$ such that $f'^* \mathcal{G}'_{x''}$ is flat over $\Ow_{X'',x''}$ is an open subset $U_{f'}$ of $X''$ by (\cite{Stacks}, 0399(2)). The $\Ow_{\widetilde{X}}$-module $\pi'^* \mathcal{G}' \cong (\pi^{*} \mathcal{F})^{\pur}$ is flat by Theorem \ref{chap1RG2}, hence by (\cite{Stacks} 02JO(3)), any point of $\widetilde{X} $ is in the inverse image of $U_{f'}$ for some $f'$. Since $\widetilde{X}$ is quasi-compact by Theorem \ref{chap1comparisoneclate}, we can find $f'$ as above such that the inverse image of $U_{f'}$ in $\widetilde{X} $ is $\widetilde{X} $. By (\cite{Stacks}, 0A2W(1)) we can assume up to taking a further refinement of $f'$ that we have $U_{f'} = X''$. Since $f'^* \mathcal{G}'$ is flat over $\Ow_{X''}$, it must coincide with its purification $(f'^{*} f^* \mathcal{F})^{\pur}$. In particular the strict transform $(f'^{*} f^* \mathcal{F})^{\pur}$ of $\mathcal{F}$ on $X''$ is flat over $\Ow_{X''}$ and finitely presented over $f'^* f^* \mathcal{A}$. This concludes the proof of \ref{chap1RGtheo2}.

\subsection{} We now prove Raynaud-Gruson's theorem \ref{chap1RGtheo}. Let $g : Y \rightarrow X, U$ and $\mathcal{F}$ be as in the statement of the theorem. Let $(X_{\lambda})_{\lambda \in \Lambda}$ and $(Y_{\lambda})_{\lambda \in \Lambda}$ be covers of $X$ and $Y$ by finitely many affine open subsets such that for each $\lambda \in \Lambda$, we have $g(Y_{\lambda}) \subseteq X_{\lambda}$. Then for each $\lambda \in \Lambda$, Theorem \ref{chap1RGtheo2} yields a finitely generated quasi-coherent ideal sheaf $\mathcal{J}_{\lambda} \subseteq \Ow_{X_{\lambda}}$ such that $V(\mathcal{J}_{\lambda})$ is disjoint from $U \cap X_{\lambda}$ and such that the blow-up $\mathcal{J}_{\lambda}$ has the desired flattening property with respect to the affine morphism $Y_{\lambda} \rightarrow X_{\lambda}$ and $\mathcal{F}_{|Y_{\lambda}}$. 

Since $\mathcal{J}_{\lambda} \Ow_{U \cap X_{\lambda}} = \Ow_{U \cap X_{\lambda}}$, there exists a unique quasi-coherent ideal sheaf $\mathcal{J}_{\lambda}'$ on $U \cup X_{\lambda}$ such that $\mathcal{J}_{\lambda |U}' = \Ow_U$ and $\mathcal{J}_{\lambda |X_{\lambda}}' = \mathcal{J}_{\lambda}$. Since $X$ and $U \cup X_{\lambda}$ are quasi-compact, there exists by (\cite{EGA1}, $9.4.7$) a quasi-coherent ideal sheaf $\mathcal{I}_{\lambda}$ of finite type on $X$ such that $\mathcal{I}_{\lambda|U \cup X_{\lambda}} = \mathcal{J}_{\lambda}'$. The blow-up of $\prod_{\lambda \in \Lambda}\mathcal{I}_{\lambda}$ then satisfies the conclusion of Theorem \ref{chap1RGtheo}.

\bibliographystyle{amsalpha}

\end{document}